\documentclass[10pt]{amsart}
\usepackage{amssymb}

\newtheorem{lemma}{Lemma}[section]

\newtheorem{theorem}[lemma]{Theorem}
\newtheorem{remark}[lemma]{Remark}
\newtheorem{proposition}[lemma]{Proposition}
\newtheorem{corollary}[lemma]{Corollary}
\newtheorem{definition}[lemma]{Definition}

\begin{document}

\title{Lipschitz-Killing curvatures and polar images}
\author{Nicolas Dutertre }
\address{Aix-Marseille Universit\'e, CNRS, Centrale Marseille, I2M, UMR 7373,
13453 Marseille, France.}
\email{nicolas.dutertre@univ-amu.fr}


\begin{abstract} 
We relate the Lipschitz-Killing measures of a definable set $X \subset \mathbb{R}^n$ in an o-minimal structure to the volumes of generic polar images. For smooth submanifolds of $\mathbb{R}^n$, such results were established by Langevin and Shifrin \cite{LangevinShifrin}.
Then we give infinitesimal versions of these results. As a corollary, we obtain a relation between the polar invariants of Comte and Merle \cite{ComteMerle} and the densities of generic polar images.
 \end{abstract}

\maketitle
\markboth{N. Dutertre}{Lipschitz-Killing curvatures and polar images}

\section{Introduction}
Let $M \subset \mathbb{R}^n$ be a smooth compact submanifold of dimension $d_M$. To each $x \in M$, we can associate a sequence of curvatures $K_0(x),\ldots,K_{d_M}(x)$ called the Lipschitz-Killing-Weyl curvatures. Let us recall their definition. We denote by $S_{N_x M}$ the unit sphere of $N_x M$, the normal space to $M$ at $x$. For $i \in \{0,\ldots,d_M\}$, the $i$-th Lipschitz-Killing-Weyl curvature $K_i(x)$ is defined by
$$K_i(x)=\int_{S_{N_x M}} \sigma_i(II_{x,v}) dv,$$
where $II_{x,v}$ is the second fundamental form on $M$ at $x$ in the direction $v$ and $\sigma_i(II_{x,v})$ is the $i$-th elementary symmetric function of its eigenvalues. The second fundamental form $II_{x,v}$ is defined on $T_x M$ as follows:
$$II_{x,v} (W_1,W_2) =-\langle \nabla_{W_1} W,W_2 \rangle,$$
for $W_1$ and $W_2$ in $T_x M$, where $\nabla$ is the covariant differentiation in $\mathbb{R}^n$ and $W$ is a local extension of $v$ normal to $M$. Note that $K_i =0$ if $i$ is odd. These curvatures are important objects: the integrals $\int_M K_i(x)dx$ appear in Weyl's tube formula \cite{Weyl} and the integral $\int_M K_{d_M}(x) dx$ in the Gauss-Bonnet formula \cite{Allendoerfer,Fenchel}. 

In \cite{LangevinShifrin}, Langevin and Shifrin explained how to compute the integrals $\int_M K_i (x) dx$, which have a differential geometric definition, by methods of differential topology. Let us explain briefly their work. Let $q$ be an even integer in $\{0,\ldots,d_M\}$, let $G_n^{d_M-q+1}$be the Grassmann manifold of $(d_M-q+1)$-dimensional linear spaces in $\mathbb{R}^n$ and for $P \in G_n^{d_M-q+1}$, let $\pi_P^M: M \rightarrow P$ be the restriction to $M$ of the orthogonal projection on $P$. Generically the discriminant $\Delta_P^M$ of $\pi_P^M$, also called the polar image of $\pi_P^M$, is almost everywhere a $(d_M-q)$-dimensional submanifold of $P$. With each regular point $y$ in $\Delta_P^M$, we can associate an integer $\alpha(\pi_P^M,y)$, which is very roughly speaking a kind of Morse index.

Langevin and Shifrin showed that
$$\int_M K_q (x) dx = {\rm cst} \frac{1}{g_n^{d_M-q+1}} \int_{G_n^{d_M-q+1}} \int_{\Delta_P^M} \alpha(\pi_P^M,y) dy dP,$$
where $g_n^{d_M-q+1}$ is the volume of $G_n^{d_M-q+1}$ and ``cst" means a universal constant that depends only on $d_M-q$ and $n$. For $q=d_M$, this formula is exactly the exchange formula or exchange principle (already proved in \cite{Langevin}) which relates the integral $\int_M K_{d_M}(x) dx$ to the Morse indices of critical points of generic orthogonal projections onto lines. 

Our first aim is to extend Langevin and Shifrin's results to a large class of non-smooth objects, namely the class of definable sets in an o-minimal structure. Definable sets are a generalization of semi-algebraic sets and global subanalytic sets, we refer the reader to classical references \cite{DriesMiller,Dries,Coste,Loi10} for basic definitions and results on this topics. The study of the geometric properties of these objects was initiated by Fu \cite{Fu94}, who developed integral geometry for compact subanalytic sets. Using the technology of the normal cycle, he associated with every compact subanalytic set $X$ of $\mathbb{R}^n$ a sequence of curvature measures
$$\Lambda_0(X,-),\ldots,\Lambda_n(X,-),$$
called the Lipschitz-Killing measures. In \cite{BroeckerKuppe} (see also \cite{BernigBroecker2,BernigBroecker}), Br\"ocker and Kuppe gave a geometric characterization of these measures using stratified Morse theory, in the more general setting of definable sets.  

Let us describe now how we adapt the technics of Langevin and Shifrin to the singular definable case. We consider a compact definable set $X \subset \mathbb{R}^n$ and assume that it is equipped with a finite definable Whitney stratification $\{S_a\}_{a \in  A}$. Let $S$ be a stratum of $X$ of dimension $d_S<n$ and let $U$ be an open subset of $X$. As in the smooth case, for $q\in \{0,\ldots,d_S\}$ and $P$ generic in $G_n^{q+1}$, the discriminant $\Delta_P^{S \cap U}$ of the restriction to $S \cap U$ of the orthogonal projection on $P$ is almost everywhere a smooth hypersurface in $P$, and to almost all $y$ in $\Delta_P^{S \cap U}$, we can assign an index $\alpha(\pi_P^S,y)$, defined by means of stratified Morse theory. Then we set (see Definition \ref{DefPolarLength1})
$$L_q(X,S,U) ={\rm cst} \int_{G_n^{q+1}} \int_{\Delta_P^{S \cap U}} \alpha(\pi_P^S,y) dy dP.$$
For $d_S < q \le n$, we set $L_q(X,S,U)=0$. If $d_S=n$, we set $L_q(X,S,U)=0$ if $q<n$ and $L_n(X,S,U)= {\rm vol}(S \cap U)$. We define the polar lengths of $X$ (see Definition \ref{DefPolarLength2}) by
$$L_q(X,U)= \sum_{a \in A} L_q (X,S_a,U) \hbox{ for } q \in \{0,\ldots,n\}.$$
In Theorem \ref{Curvature=Length}, we establish a connection between the Lipschitz-Killing measures and the polar lengths, namely we prove that for any open subset $U$ of $X$ and for $q \in \{0,\ldots,n\}$,
$$\Lambda_q(X,U)= L_q(X,U).$$
Our second goal is to give infinitesimal versions of the previous equalities.  We consider $(X,0) \subset (\mathbb{R}^n,0)$ a germ of a closed definable set equipped with a finite definable Whitney stratification. We define localized versions of the polar lengths that we denote by $L_k^{\rm loc}(X,0)$, $k=0,\ldots,n$. Roughly speaking, these localized polar lengths are mean-values over Grassmann manifolds of weighted sums of densities of polar images. We prove that (see Theorem \ref{Curvature=LengthLocal}) 
$$\lim_{\epsilon \rightarrow 0} \frac{\Lambda(X,X \cap B_\epsilon)}{b_k \epsilon^k} = L_k^{\rm loc}(X,0),$$
for $ k \in \{0,\ldots,n\}$. As a consequence, we obtain in Corollary \ref{PolarLocalAndPolarInv} a relation between the $L_k^{\rm loc}(X,0)$'s and the polar invariants of Comte and Merle \cite{ComteMerle}. This implies that the $L_k^{\rm loc}(X,0)$'s are continuous along the strata of a Verdier stratification of $X$. This result should be related to Teissier's famous result on the constancy of polar multiplicities along the strata of a Whitney stratification of a complex analytic set \cite{Teissier}.

In complex analytic and algebraic geometry, polar and relative polar varieties and their relations with curvatures, characteristic classes and equisingularity problems have been widely studied by many authors since the 70's. We cannot give here a complete list of all the interesting papers published on these subjects and we apologize for this. The study of real polar varieties and real equisingularity, suggested by Trotman \cite{Trotman}, was started by Comte in \cite{Comte}. It was continued in \cite{ComteMerle}, \cite{Valette} and \cite{NguyenValette}. We hope that the present paper will provide new interesting developments in this theory. 

Throughout the paper, we will use the following notations and conventions (some of them have already appeared in this
introduction):
\begin{itemize}
\item $s_k$ is the volume of unit sphere $S^k$ of dimension $k$ and $b_k$ is the volume of the unit ball $B^k$ of dimension $k$,
\item $\beta(n,k)= \frac{\Gamma(\frac{k+1}{2})\Gamma(\frac{n-k+1}{2})}{\Gamma(\frac{1}{2})\Gamma(\frac{n+1}{2})}$ where $\Gamma$ is the Euler function,
\item for $k \in \{0,\ldots,n \}$, $G_{n}^k$ is the Grassmann manifold of $k$-dimension  linear spaces in $\mathbb{R}^{n}$ equipped with the $O(n)$-invariant density (see for instance \cite{Santalo}, p.200),
$g_n^k$ is its volume,
\item $A_{n}^k$ is the affine grassmanian of $k$-dimensional affine spaces in $\mathbb{R}^{n}$,
\item if $P$ is a linear subspace of $\mathbb{R}^{n}$, $S_P$ is the unit sphere in $P$, $G_P^k$  is  the Grassmann manifold of
$k$-dimensional linear spaces in $P$, $A_P^k$ is the affine grassmanian of $k$-dimensional affine spaces in $P$, $P^\perp$ is the orthogonal space to $P$, $\pi_P : \mathbb{R}^{n} \rightarrow P$ is the orthogonal projection on $P$ and for any subset $A \subset \mathbb{R}^n$, $\pi_P^A : A \rightarrow P$ is its restriction to $A$,
\item for $v \in \mathbb{R}^n$, the function $v^* : \mathbb{R}^n \rightarrow \mathbb{R}$ is defined by $v^*(y)= \langle v, y \rangle$,
\item in $\mathbb{R}^n$, $B_\epsilon(x)$ is the closed ball of radius $\epsilon$ centered at $x$ and $S_\epsilon(x)$ is the sphere of radius $\epsilon$ centered at $x$, if $x=0$, we simply write $B_\epsilon$ and $S_\epsilon$,
\item if $X \subset \mathbb{R}^{n}$,   ${\rm Sing}(X)$ is the singular set of $X$, $\overline{X}$ is its
topological closure, $\mathring{X}$ its topological interior, ${\rm Fr}(X)=\overline{X} \setminus X$ its frontier,
\item  when it makes sense, ${\rm vol}(X)$ means the volume of the set $X$ and $d_X$ its dimension,
\item if $v_1,\ldots,v_k$ are vectors in $\mathbb{R}^{n}$, $[v_1,\ldots,v_k]$ is the linear
space spanned by $v_1,\ldots,v_k$,
\item a universal constant that we do not want to specify will be denoted by ``cst",
\item the word ``smooth" means of class $C^3$ at least,
\item if $f(x,y)=f(x_1,\ldots,x_n,y_1,\ldots,y_m)$ is a smooth function, $\nabla f$ is its gradient and $\nabla_x f$ is its gradient with respect to the variables $x_1,\ldots,x_n$.
\end{itemize}

The paper is organized as follows. Section 2 contains a summary of results on stratified mappings, stratified More theory and the definitions of the Lipschitz-Killing measures and the polar invariants. In Section 3, we study generic polar varieties and images and prove the relation between the Lipschitz-Killing measures and the polar lengths. Section 4 deals with the local situation and contains infinitesimal versions of the results of Section 3.

The author thanks David Mond and Terry Gaffney for interesting discussions on double points of fold singularities.

\section{Stratified mappings, stratified Morse theory and Lipschitz-Killing curvatures}
In this section, we present different mathematical objects, notions and results that we will use in the next sections.

\subsection{Stratified mappings}
In this subsection, we recall well-known facts on mappings defined on stratified sets and on their critical points and values.

Let $X \subset \mathbb{R}^n$ be a compact definable set equipped with a finite definable Whitney stratification $\mathcal{S}=\{ S_a \}_{a \in A}$.  
The fact that such a stratification exists is due to Loi \cite{Loi98}. Recently Nguyen, Trivedi and Trotman \cite{NguyenTrivediTrotman} gave another proof of this result. Let $f : \mathbb{R}^n \rightarrow \mathbb{R}^k$ be a smooth definable mapping and let $f_{\vert X} : X \rightarrow \mathbb{R}^k$ be its restriction to $X$. We assume further that $f$ is a submersion in an open neighborhood of $X$.

A point $x$ in $X$ is a (stratified) critical point of $f_{\vert X}$ if $x$ is a critical point of $f_{\vert S}$, where $S$ is the stratum that contains $x$. This means that $${\rm dim}(T_x S + T_x f^{-1} (f(x))) <n,$$ or equivalently $${\rm dim}(T_x S \cap T_x f^{-1} (f(x))) > d_S-k.$$ Note that if $d_S <k$, then all the points of $S$ are critical. 
If $x$ is not a critical point of $f_{\vert X}$, we say that $x$ is a regular point of $f_{\vert X}$. We denote by $\Sigma_f^S$ the set of critical points of $f_{\vert S}$ and we set $\Sigma_f^X = \cup_{a \in A} \Sigma_f^{S_a}$. 
\begin{lemma}\label{Transverse0}
The set $\Sigma_f^X$ is a compact definable set of $X$.
\end{lemma}
\proof It is enough to prove that it is closed. Let $(x_m)_{m \in \mathbb{N}}$ be a sequence of points in $\Sigma_f^X$ that tends to $x$. We can assume that $(x_m)_{m \in \mathbb{N}}$ is included in $\Sigma_f^S$. If $x$ belongs to $S$ then $x$ belongs to $\Sigma_f^S \subset \Sigma_f^X$. If $x$ belongs to $S' \subset {\rm Fr}(S)$ then, by Whitney condition (a), there exists a vector space $T$ of dimension $d_S$ such that ${\rm dim}(T + T_x f^{-1} (f(x))) <n$ and such that $T_x S' \subset T$. Therefore $x \in \Sigma_f^{S'}$. \endproof

Let $y \in \mathbb{R}^k$. We say that $y$ is a critical value of $f_{\vert X}$ if $f_{\vert X}^{-1}(y)$ contains a critical point of $f_{\vert X}$. Otherwise we say that $y$ is a regular value of $f_{\vert X}$. Note that if $y$ is a regular value of $f_{\vert X}$, then $f^{-1}(y) \cap X$ is a Whitney stratified set (see for instance \cite{OrroTrotman}). 
\begin{lemma}\label{Transverse}
The set of regular values of $f_{\vert X}$ is an  open definable and dense subset of $\mathbb{R}^k$.
\end{lemma}
\proof This set is equal to $\mathbb{R}^k \setminus f(\Sigma_f^X)$. By Bertini-Sard's theorem (see \cite{BochnakCosteRoy}), each set $f(\Sigma_f^S)$ is a definable set of dimension less or equal to $k-1$. Hence $f(\Sigma_f^X)$ is a compact definable set of dimension less or equal to $k-1$. \endproof 
\subsection{Stratified Morse theory}
Let $M$ be a smooth riemannian manifold.
Let $X$ be a Whitney-stratified set of $M$. Let $x$ be a point in $X$ and let $S(x)$ be the stratum that contains $x$. A generalized tangent space at $x$ is a
limit of a sequence of tangent spaces $(T_{y_m} S_1)_{m \in \mathbb{N}}$, where $S_1$ is a stratum distinct from $S(x)$ such that $x \in
\overline{S_1}$ and $(y_m)_{m \in \mathbb{N}}$ is a sequence of points in $S_1$ tending to $x$. 

A Morse function $f :X \rightarrow \mathbb{R}$ is the restriction of a smooth function $\tilde{f} : M \rightarrow
\mathbb{R}$ such that the following conditions hold:
\begin{enumerate}
\item $f$ is proper and the critical values of $f$ are distinct i.e., if $x$ and $y$ are two distinct critical points of $f$ then $f(x) \not=
f(y)$,
\item for each stratum $S$ of $X$, the critical points of $f_{\vert S}$ are non-degenerate,
\item at each critical point $x$, the differential $D\tilde{f}(x)$ does not annihilate any generalized tangent
space at $x$.
\end{enumerate} 

Let $f : X \rightarrow \mathbb{R}$ be a Morse function, $S$ a stratum of $X$ and $x$ a critical point of $f_{\vert S}$. Let us write
$c=f(x)$.  The local Morse data for $f$ at $x$ is the pair 
$$\left( B^M_\epsilon(x) \cap f^{-1}([c-\delta,c+\delta]), B^M_\epsilon(x) \cap f^{-1}(c-\delta) \right),$$
where $0 < \delta \ll \epsilon \ll 1$ and $B^M_\epsilon(x)$ is the $d_M$-dimensional closed ball centered at $x$ and of radius $\epsilon$ in $M$. This definition is justified by the following property (see \cite{GoreskyMacPherson} or \cite{Hamm}). There exists an
open subset $A$ in $\mathbb{R}^+ \times \mathbb{R}^+$ such that 
\begin{enumerate}
\item the closure $\overline{A}$ of $A$ in $\mathbb{R}^2$ contains an interval $]0,\epsilon[$, $\epsilon >0$, such that for all $a \in
]0,\epsilon[$, the set $\{ b \in \mathbb{R}^+ \  \vert \ (a,b) \in A \}$ contains an open interval $]0,\delta(a)[$ with
$\delta(a)>0$,
\item for all $(\epsilon,\delta) \in A$, the above pairs of spaces are homeomorphic.
\end{enumerate}
The local Morse data are Morse data in the sense that $f^{-1}(]-\infty, c + \delta])$ is homeomorphic to the space one gets by
attaching $B_\epsilon^M(x) \cap f^{-1}([c-\delta,c+\delta])$ at $f^{-1}(]-\infty, c - \delta])$ along $B_\epsilon^M(x) \cap
f^{-1}(c-\delta)$ (see \cite{GoreskyMacPherson}, I 3.5). 

If $x$ belongs to a stratum of dimension $d_X$, the local Morse data at $x$ are homeomorphic to the classical Morse data 
$(B^\lambda \times B^{d_X-\lambda}, \partial B^\lambda \times B^{d_X-\lambda})$ 
where $\lambda$ is the Morse index of $f$ at $x$.  

If $x$ belongs to a zero-dimensional stratum then $B_\epsilon^M(x) \cap f^{-1}([c-\delta,c+\delta])$ has the
structure of a cone (see \cite{GoreskyMacPherson}, I 3.11). 

If $x$ lies in a stratum $S$ with $0< d_S < d_X$, then one can consider
the classical Morse data of $f_{\vert S}$ at $x$. We will call them tangential Morse data and denote them by
$(P_{\rm tg},Q_{\rm tg})$. One may choose a normal slice of $N$ at $x$, that is a closed submanifold of $M$ of
dimension $d_M-d_S$, which intersects $S$ in $x$ transversally. We define the normal Morse data
$(P_{\rm nor},Q_{\rm nor})$ at $x$ to be the local Morse data of $f_{\vert X \cap N}$ at $x$. Goresky and Mac-Pherson proved that 
\begin{itemize}
\item the normal Morse data are well defined, that is to say they are independent of the Riemannian metric and the choice of
the normal slice (\cite{GoreskyMacPherson}, I 3.6),
\item the local Morse data $(P,Q)$ of $f$ at $x$ are the product of the tangential and the normal Morse data (\cite{GoreskyMacPherson}, I 3.7):
$$(P,Q) \simeq (P_{\rm tg} \times P_{\rm nor}, P_{\rm tg} \times Q_{\rm nor} \cup P_{\rm nor} \times Q_{\rm  tg}).$$
\end{itemize}
This implies that $P=B_\epsilon^M(x) \cap f^{-1}([c-\delta,c+\delta])$ has the structure of a cone. 

\begin{definition}
{\rm Let $x \in X$ be a critical point of the Morse function $f : X \rightarrow \mathbb{R}$. Let $(P,Q)$ be the local Morse data of $f$ at
$x$. The Euler-Poincar\'e characteristic $\chi(P,Q)=1-\chi(Q)$ is called the stratified Morse index of $f$ at $x$ and is denoted by
${\rm ind}(f,X,x)$. If $x \in X$ is not a critical point of $f$, we set ${\rm ind}(f,X,x)=0$.}
\end{definition}

If $(P_{tg},Q_{tg})$ and $(P_{nor},Q_{nor})$ are the tangential and normal Morse data, then one has
$$\chi(P,Q)=\chi(P_{tg},Q_{tg}) \cdot \chi(P_{nor},Q_{nor}).$$
Thus we can write $${\rm ind}(f,X,x)={\rm ind}_{\rm tg}(f,X,x) \cdot {\rm ind}_{\rm nor}(f,X,x),$$ where ${\rm ind}_{\rm tg}(f,X,x)$ is called the tangential
Morse index and ${\rm ind}_{\rm nor}(f,X,x)$ the normal Morse index. 
Note that if $x$ belongs to a zero-dimensional stratum then ${\rm ind}_{\rm tg} (f,X,x)=1$. If $x$ belongs a stratum of dimension $d_X$, then ${\rm ind}_{\rm tg}(f,X,x)=(-1)^\lambda$ and ${\rm ind}_{\rm nor}(f,X,x)=1$, where $\lambda$ is the Morse index of $f$ at $x$.

The following theorem relates the Euler characteristic of $X$ to the indices of the critical points of a Morse function.
\begin{theorem}
Let $X \subset M$ be a compact Whitney-stratified set and let $f : X \rightarrow \mathbb{R}$ be a Morse function. 
One has 
$$\chi(X)=\sum_{x \in X} {\rm ind}(f,X,x).$$
\end{theorem}

\subsection{Lipschitz-Killing measures of definable sets}
In this subsection, we present the Lipschitz-Killing measures of a definable set in an o-minimal structure.  We describe Br\"ocker and Kuppe's approach \cite{BroeckerKuppe}.

Let $X \subset \mathbb{R}^n$ be a compact definable set equipped with a finite definable Whitney stratification $\mathcal{S}=\{ S_a \}_{a \in A}$.  

Let  us fix a stratum $S$.
For $k \in \{0,\ldots,d_S\}$, let $\lambda_k^S : S \rightarrow \mathbb{R}$ be defined by
$$ \lambda_k^S(x) = \frac{1}{s_{n-k-1}} \int_{S_{T_x S^\perp} }{\rm ind}_{\rm nor}(v^*,X,x) \sigma_{d_S-k} (II_{x,v}) dv,$$
where $II_{x,v}$ is the second fundamental form on $S$ in the direction of $v$ and where $\sigma_{d_S-k} (II_{x,v})$ is the $(d_S-k)$-th elementary symmetric function of its eigenvalues. The index ${\rm ind}_{\rm nor}(v^*,X,x)$ is defined as follows:
$${\rm ind}_{\rm nor}(v^*,X,x)= 1-\chi \Big( X \cap N_x \cap B_\epsilon(x) \cap \{ v^*= v^*(x)-\delta \} \Big),$$
where $0 < \delta \ll \epsilon \ll 1$ and $N_x$ is a normal (definable) slice to $S$ at $x$ in $\mathbb{R}^n$ . 
Since we work in the definable setting, this index is well-defined thanks to Hardt's theorem \cite{Hardt,Coste}. Furthermore when $v^*_{\vert X}$ has a stratified Morse critical point at $x$, it coincides with the normal Morse index at $x$ of a function $f : \mathbb{R}^n \rightarrow \mathbb{R}$ such that $f_{\vert X}$ has a stratified Morse critical point at $x$ and $\nabla f (x) = v$. For $k \in \{d_S+1,\ldots,n\}$, we set $\lambda_k^S(x)=0$. 

If $S$ has dimension $n$ then for all $x \in S$, we put $\lambda_0^S(x)=\cdots=\lambda_{n-1}^S(x)=0$ and $\lambda_n^S(x)=1$. If $S$ has dimension $0$ then we set
$$\lambda_0^S(x)= \frac{1}{s_{n-1}} \int_{S^{n-1}}  {\rm ind}_{\rm nor}(v^*,X,x) dv,$$
and $\lambda_k^S(x)=0$ if $k>0$. 
\begin{definition}
{\rm For every Borel set $U \subset X$ and for every $k \in \{0,\ldots,n\}$, we define $\Lambda_k(X,U)$ by
$$\Lambda_k(X,U)= \sum_{a \in A} \int_{S_a \cap U} \lambda_k^{S_a} (x) dx.$$}
\end{definition}
These measures $\Lambda_k(X,-)$ are called the Lipschitz-Killing measures of $X$. Note that for any Borel set $U$ of $X$, we have
$$\Lambda_{d_X+1}(X,U)= \cdots=\Lambda_n(X,U)=0,$$ 
and $\Lambda_{d_X}(X,U)= \mathcal{L}_{d_X}(U)$, where $\mathcal{L}_{d_X}$ is the $d_X$-th dimensional Lebesgue measure in $\mathbb{R}^n$. 
If $X$ is smooth then for $k \in \{0,\ldots,d_X \}$, $\Lambda_k(X,U)$ is equal to 
$$\frac{1}{s_{n-k-1} } \int_U K_{d_X-k}(x) dx.$$

As in the smooth case, the measure $\Lambda_0(X,-)$ satisfies an exchange formula (see \cite{BroeckerKuppe}).
\begin{proposition}\label{ExchangeFormula}
For every Borel set $U \subset X$, we have 
$$\Lambda_0(X,U) =\frac{1}{s_{n-1}} \int_{S^{n-1}} \sum_{x \in X}  {\rm ind} (v^*,X,x) dv.$$
\end{proposition}
In the next section, we will generalize this exchange formula to the other measures $\Lambda_k(X,-)$, $k \ge 1$. 

The Lipschitz-Killing measures satisfy the kinematic formula (see \cite{Fu94,BroeckerKuppe,BernigBroecker2}). We will need a particular case of this formula, namely the linear kinematic formula.
\begin{proposition}\label{GlobalKinematicFormula}
Let $U$ be an open subset of $X$. For $k \in \{0,\ldots,n \}$, we have
$$\Lambda_{n-k}(X,U) = {\rm cst} \int_{A_n^k} \Lambda_0(X \cap E, X \cap E \cap U) dE.$$
\end{proposition}
\proof See \cite{BroeckerKuppe}, Corollary 8.5. \endproof

In \cite{DutertreProcTrotman}, we gave a localized  version of this equality. Let $(X,0) \subset (\mathbb{R}^n,0)$ be the germ of a closed definable set. Let $H \in G_n^{n-k}$, $k \in \{1,\ldots, n \}$, and let $v$ be an element in $S_{H^\perp}$. For $\delta >0$, we denote by $H_{v,\delta}$ the $(n-k)$-dimensional affine space $H+\delta v$ and we set 
$$\beta_0(H,v) = \lim_{\epsilon \rightarrow 0} \lim_{\delta \rightarrow 0} \Lambda_0(H_{\delta,v} \cap X, H_{\delta,v} \cap X \cap B_\epsilon).$$
Then we set 
$$\beta_0 (H) =\frac{1}{s_{k-1}}\int_{S_{H^\perp}} \beta_0 (H,v) dv.$$
\begin{theorem}\label{LocalKinematicFormula}
For $k \in \{1,\ldots,n \}$, we have 
$$ \lim_{\epsilon \rightarrow 0} \frac{\Lambda_k(X,X \cap B_\epsilon)}{b_k \epsilon^k } = \frac{1}{g_n^{n-k}} \int_{G_n^{n-k}} \beta_0 (H) dH.$$
\end{theorem}
\proof See \cite{DutertreProcTrotman}, Theorem 5.5. \endproof

In \cite{DutertreProcTrotman}, we also established a relation between the limits $\lim_{\epsilon \rightarrow 0} \frac{\Lambda_k(X,X \cap B_\epsilon)}{b_k \epsilon^k }$ and the polar invariants introduced by Comte and Merle in \cite{ComteMerle}.  These polar invariants are real versions of the vanishing Euler characteristics of L\^e and Teissier (see \cite{LeTeissierAnnals}, Proposition 6.1.8 or \cite{LeTeissierArcata}, (3.1.2)). They can be defined as follows. Let $H \in G_n^{n-k}$, $k \in \{1,\ldots, n \}$, and let $v$ be an element in $S_{H^\perp}$. We set 
$$\alpha_0(H,v) = \lim_{\epsilon \rightarrow 0} \lim_{\delta \rightarrow 0} \chi(H_{\delta,v} \cap X \cap B_\epsilon),$$
$$\alpha_0(H)=\frac{1}{s_{k-1} }\int_{S_{H^\perp}} \alpha_0 (H,v) dv,$$
and then 
$$\sigma_k (X,0)=\frac{1}{g_n^{n-k}} \int_{G_n^{n-k}} \alpha_0 (H) dH.$$
Moreover, we put $\sigma_0 (X,0)=1$. 
\begin{theorem}\label{CurvAndPolar}
For $k \in \{0,\ldots,n-1\}$, we have
$$\lim_{\epsilon \rightarrow 0} \frac{\Lambda_k(X,X \cap B_\epsilon)}{b_k \epsilon^k} = \sigma_k(X,0) -\sigma_{k+1}(X,0).$$
Furthermore, we have
$$\lim_{\epsilon \rightarrow 0} \frac{\Lambda_n(X,X \cap B_\epsilon)}{b_n \epsilon^n} = \sigma_n(X,0).$$
\end{theorem}
\proof See \cite{DutertreProcTrotman}, Theorem 5.6. \endproof

\section{Curvature measures and volume of polar images}\label{Global}
In this section,  we relate the curvature measures of definable sets to the volumes of polar images of generic projections. 
We start recalling the definition of polar varieties and polar images.

\begin{definition}\label{defpolar}
{\rm Let $P \in G_n^k$, $k=1,\ldots,n$, and let $M \subset \mathbb{R}^n$ be a smooth submanifold.
The polar variety (or polar set) $\Sigma_P^M$ is the set of critical points of $\pi_P^M$, i.e.,
$$\Sigma_P^M = \left\{ x \in M \ \vert \ {\rm dim}(T_x M \cap P^\perp) \ge d_M-k +1 \right\},$$
if $k \le d_M$. If $k > d_M$, we set $\Sigma_P^M = M$.

The polar image of $\pi_P^M$ is the set $\Delta_P^M$ defined by $\Delta_P^M = \pi_P^M(\Sigma_P^M)$.}
\end{definition}

Let $X \subset \mathbb{R}^n$ be a compact definable set in an o-minimal structure. We equip it with a definable  Whitney stratification $\mathcal{S}$. The following statements describe the structure of $\Sigma_P^S$ and $\Delta_P^S$ for a stratum $S$ of $X$. 

The next two lemmas can be proved using the machinery of modular submanifolds of multi-jet spaces developed by Mather and the fact the Thom-Boardman manifolds are modular (see \cite{Mather}, Theorem 2). The proofs we present below are more elementary and enable us to stay in the definable setting.
\begin{lemma}\label{GenericProj}
Let $S$ be a stratum of $X$ such that $d_S< n$. If $k \le d_S+1$ then for almost all $P \in G_n^k$, $\Sigma_P^S$ is either empty or a $(k-1)$-dimensional definable set. Moreover the set 
$$\Sigma_P^{'S}= \left\{ x \in \Sigma_P^S \ \vert \ {\rm dim}(T_x S \cap P^\perp) \ge d_S-k +2 \right\},$$
is a definable subset of $\Sigma_P^S$ of dimension at most $k-2$ and $\Sigma_P^S \setminus \Sigma_P^{'S}$ is smooth.
\end{lemma}
\proof We treat first the case $k \le d_S$. Let $\mathcal{V} \subset (\mathbb{R}^n)^k$ be the open subset consisting of the $k$-tuples $(v_1,\ldots,v_k)$ such that $(v_1,\ldots,v_k)$ has rank $k$. Let $M$ be the following definable set:

$$\displaylines{
M = \Big\{ (x,v_1,\ldots,v_k) \in \mathbb{R}^n \times \mathcal{V} \ \vert \  x \in S \hbox{ and dim}(T_x S \cap P^\perp) \ge d_s-k+1  \hfill \cr
 \hfill    \hbox{ where } P=[ v_1,\ldots, v_k ] \Big\}.  \qquad \cr
 }$$
The set $M$ is the finite disjoint union of the sets $M_i$, $i=1,\ldots,k$, where 
$$\displaylines{
  M_i = \Big\{ (x,v_1,\ldots,v_k) \in \mathbb{R}^n \times \mathcal{V} \ \vert \  x \in S \hbox{ and dim} (T_x S \cap P^\perp) = d_s-k+i    \hfill \cr
\hfill   \hbox{ where } P=[ v_1,\ldots, v_k ] \Big\}. \qquad \cr
}$$
The set $M_1$ is a definable manifold of dimension $nk+k-1$. To see this, we can assume that locally $S$ is given by
$$S = \{ x \in \mathbb{R}^n \ \vert \ f_1(x)= \cdots=f_{c_S}(x)= 0 \},$$
where $c_S=n-d_S$ and $f_1,\ldots, f_{c_S}$ are smooth definable functions. From now on, we set $V=(v_1,\ldots,v_k)$. Then $(x,V)$ is in $M_1$ if and only if the matrix
$$\left( \begin{array}{ccc}
\frac{\partial f_1}{\partial x_1} (x) & \cdots & \frac{\partial f_1}{\partial x_n} (x)\cr
\vdots & \ddots & \vdots \cr
\frac{\partial f_{c_S}}{\partial x_1} (x)& \cdots & \frac{\partial f_{c_S}}{\partial x_n} (x)\cr
v_1^1 & \cdots & v_1^n \cr
\vdots & \ddots & \vdots \cr
v_k^1 & \cdots & v_k^n \cr  
\end{array} \right),$$
has rank $c_S+k-1$, where we use the notation $v_i=(v_i^1, \ldots, v_i^n) \in \mathbb{R}^n$.  It is an easy exercise of linear algebra to see that we can assume that the matrix 
$$\left( \begin{array}{ccc}
\frac{\partial f_1}{\partial x_1}(x) & \cdots & \frac{\partial f_1}{\partial x_n}(x) \cr
\vdots & \ddots & \vdots \cr
\frac{\partial f_{c_S}}{\partial x_1}(x) & \cdots & \frac{\partial f_{c_S}}{\partial x_n}(x) \cr
v_1^1 & \cdots & v_1^n \cr
\vdots & \ddots & \vdots \cr
v_{k-1}^1 & \cdots & v_{k-1}^n \cr  
\end{array} \right),$$
has rank $c_S+k-1$ and so that the minor 
$$\left\vert \begin{array}{ccc}
\frac{\partial f_1}{\partial x_1} (x) & \cdots & \frac{\partial f_1}{\partial x_{c_S+k-1}} (x)\cr
\vdots & \ddots & \vdots \cr
\frac{\partial f_{c_S}}{\partial x_1}(x) & \cdots & \frac{\partial f_{c_S}}{\partial x_{c_S+k-1}} (x) \cr
v_1^1 & \cdots & v_1^{c_S+k-1} \cr
\vdots & \ddots & \vdots \cr
v_{k-1}^1 & \cdots & v_{k-1}^{c_S+k-1} \cr  
\end{array} \right\vert,$$
does not vanish. Therefore the matrix 
$$\left( \begin{array}{ccc}
\frac{\partial f_1}{\partial x_1} (x)& \cdots & \frac{\partial f_1}{\partial x_n}(x) \cr
\vdots & \ddots & \vdots \cr
\frac{\partial f_{c_S}}{\partial x_1}(x) & \cdots & \frac{\partial f_{c_S}}{\partial x_n} (x)\cr
v_1^1 & \cdots & v_1^n \cr
\vdots & \ddots & \vdots \cr
v_{k-1}^1 & \cdots & v_{k-1}^{n}  \cr  
\end{array} \right),$$
has rank $c_S+k-1$ if and only if the following minors:
$$m_j =\left\vert \begin{array}{cccc}
\frac{\partial f_1}{\partial x_1}(x) & \cdots & \frac{\partial f_1}{\partial x_{c_S+k-1}} (x)& \frac{\partial f_1}{\partial x_j} (x)\cr
\vdots & \ddots & \vdots \ & \vdots \cr
\frac{\partial f_{c_S}}{\partial x_1}(x) & \cdots & \frac{\partial f_{c_S}}{\partial x_{c_S+k-1}}  (x)& \frac{\partial f_{c_S}}{\partial x_j}v \cr
v_1^1 & \cdots & v_1^{c_S+k-1}  & v_1^j\cr
\vdots & \ddots & \vdots & \vdots \cr
v_{k}^1 & \cdots & v_{k}^{c_S+k-1} & v_k^j  \cr  
\end{array} \right\vert,$$
$j = c_S+k,\ldots, n$, vanish. Hence $M_1$ is defined locally by $c_S+d_S-k+1=n-k+1$ equations. Looking at the partial derivatives of the $m_j$'s with respect to the variables $v_k^j$, $j = c_S+k,\ldots, n$, it easy to see that the gradient vectors of the functions defining $M_1$ are linearly independent, because the above minor of size $c_S+k-1$ is not zero. Hence $M_1$ is a submanifold of dimension $nk+k-1$.  Now let us consider the projection $\pi_1 : M_1 \rightarrow \mathcal{V}, (x,V) \mapsto V$. Bertini-Sard's theorem implies that the set $D_{\pi_1}$ of critical values of $\pi_1$ is a definable set of dimension  less than $nk$. Hence for all $V \in \mathcal{V} \setminus D_{\pi_1}$, $\pi_1^{-1}(V)$ is a smooth definable set of dimension $k-1$ (possibly empty). But $\pi_1^{-1}(V)$ is exactly the set 
$$\left\{ x \in \Sigma_P^S \ \vert \ {\rm dim}(T_x S \cap P^\perp) = d_S-k +1 \right\},$$
where $V=(v_1,\ldots,v_k)$ and $P=[ v_1,\ldots,v_k ]$. 
Similarly we can show that the sets $M_i$, $i \ge 2$, are definable manifolds of dimension strictly less than $nk+k-1$ and we can consider the projections $\pi_i : M_i \rightarrow \mathcal{V}, (x,V) \mapsto V$. As above, for all $V \in \mathcal{V} \setminus D_{\pi_i}$, where $D_{\pi_i}$ is the discriminant of $\pi_i$, the set 
$$\left\{ x \in \Sigma_P^S \ \vert \ {\rm dim}(T_x S \cap P^\perp) = d_S-k +i \right\},$$
is a smooth definable set of dimension less or equal to $k-2$. Let us denote by ${\bf D}$ the union of the discriminants $D_{\pi_i}$. Then for all $V$ in $\mathcal{V} \setminus {\bf D}$, $\Sigma_P^S$ is  a $(k-1)$-dimensional definable set,
$\Sigma_P^{'S}$
is a definable subset of $\Sigma_P^S$ of dimension at most $k-2$ and $\Sigma_P^S \setminus \Sigma_P^{'S}$ is smooth, where $V=(v_1,\ldots,v_k)$ and $P=[v_1,\ldots,v_k ]$. 

The same proof works if $k=d_S+1$. \endproof
When $k \le d_S$, let $\Sigma_P^{S,0}$ be set of fold points of $\pi_P^S$ and when $k=d_S+1$, let $\Sigma_P^{S,0}$ be the set of regular points of $\pi_P^S$. The following lemma gives a more precise description of $\Sigma_P^S$.
\begin{lemma}\label{GenericFolds}
Let $S$ be a stratum of $X$  such that $d_S< n$. If $k \le d_S+1$ then for almost all $P \in G_n^k$, $\Sigma_P^{S,0}$ is a smooth definable set of dimension $k-1$ and $\Sigma_P^S$  admits the following decomposition:
$$\Sigma_P^S = \Sigma_P^{S,0} \sqcup Z_P^S,$$
where $Z_P^S$ is a closed definable  set of dimension strictly less than $k-1$.
\end{lemma}
\proof Let us treat first the case $k \le d_S$. We recall that $x$ is a fold point of $\pi_P^S$ if ${\rm dim}(T_x S \cap P^\perp) = d_S-k+1$ and ${\rm dim}(T_x \Sigma_P^S \cap P^\perp)=0$. We keep the notations of the previous lemma. Let $N$ be the following definable set:
$$\displaylines{
N = \Big\{ (x,v_1,\ldots,v_k) \in \mathbb{R}^n \times( \mathcal{V} \setminus \overline{ {\bf D}})  \ \vert \  x \in S, \ \hbox{dim}(T_x S \cap P^\perp) = d_s-k+1  \hfill \cr
 \hfill  \hbox{and dim} (T_x \Sigma_P^S \cap P^\perp) \ge 1  \hbox{ where } P=[ v_1,\ldots, v_k ] \Big\}.  \qquad \cr
 }$$
The set $N$ is the finite union of the sets $N_i$, $i=1,\ldots, k-1$, where 
$$\displaylines{
N_i = \Big\{ (x,v_1,\ldots,v_k) \in \mathbb{R}^n \times( \mathcal{V} \setminus \overline{ {\bf D}})  \ \vert \  x \in S, \ \hbox{dim}(T_x S \cap P^\perp) = d_s-k+1  \hfill \cr
 \hfill  \hbox{and dim} (T_x \Sigma_P^S \cap P^\perp)= i  \hbox{ where } P=[v_1,\ldots, v_k ] \Big\}.  \qquad \cr
 }$$
The set $N_1$ is a definable submanifold of dimension $nk+k-2$. To see this, let $(x,V)$ be a point  in $M_1$. Since ${\rm dim}(T_x S \cap P^\perp)=d_S-k+1$, by the proof of Lemma \ref{GenericProj} we can assume that around $(x,V)$, $M_1$ is defined by the vanishing of smooth definable functions 
$f_1,\ldots,f_{c_S}$ and $m_{c_S+k},\ldots,m_n$, where $f_1,\ldots,f_{c_S}$ depend on $x$, $m_{c_S+k},\ldots,m_n$ depend on $x$ and $V$ and where 
$${\rm rank}(\nabla f_1(x),\ldots,\nabla f_{c_S}(x),v_1,\ldots,v_{k-1})=c_S+k-1.$$
Furthermore since $V \notin D_{\pi_1}$, $\Sigma_P^S$ is defined locally at $x$ by the vanishing of $f_1,\ldots,f_{c_S}$ and $m_{c_S+k}(-,V),\ldots, m_n(-,V)$ and the following gradient vectors:
$$\nabla f_1,\ldots,\nabla f_{c_S}, \nabla_x m_{c_S+k}(-,V), \ldots, \nabla_x m_n (-,V),$$
are linearly independent.  So we see that $(x,V)$ belongs to $N_1$ if and only if the matrix
$$\left( \begin{array}{ccc}
\frac{\partial f_1}{\partial x_1} (x) & \cdots & \frac{\partial f_1}{\partial x_n} (x)\cr
\vdots & \ddots & \vdots \cr
\frac{\partial f_{c_S}}{\partial x_1} (x) & \cdots & \frac{\partial f_{c_S}}{\partial x_n} (x)\cr
\frac{\partial m_{c_S+k}}{\partial x_1} (x)& \cdots & \frac{\partial m_{c_S+k}}{\partial x_n} (x)\cr
\vdots & \ddots & \vdots \cr
\frac{\partial m_n}{\partial x_1} (x) & \cdots & \frac{\partial m_n}{\partial x_n} (x)\cr
v_1^1 & \cdots & v_1^n \cr
\vdots & \ddots & \vdots \cr
v_k^1 & \cdots & v_k^n \cr  
\end{array} \right),$$
has rank $n-1$. Since ${\rm rank}(\nabla f_1(x),\ldots,\nabla f_{c_S}(x), v_1,\ldots,v_{k-1} )= c_S+k-1$, this condition is equivalent to the fact that the matrix 
$$C= \left( \begin{array}{ccc}
\frac{\partial f_1}{\partial x_1} (x) & \cdots & \frac{\partial f_1}{\partial x_n} (x)\cr
\vdots & \ddots & \vdots \cr
\frac{\partial f_{c_S}}{\partial x_1} (x) & \cdots & \frac{\partial f_{c_S}}{\partial x_n} (x)\cr
\frac{\partial m_{c_S+k}}{\partial x_1} (x)& \cdots & \frac{\partial m_{c_S+k}}{\partial x_n} (x)\cr
\vdots & \ddots & \vdots \cr
\frac{\partial m_n}{\partial x_1} (x) & \cdots & \frac{\partial m_n}{\partial x_n} (x)\cr
v_1^1 & \cdots & v_1^n \cr
\vdots & \ddots & \vdots \cr
v_{k-1}^1 & \cdots & v_{k-1}^n \cr 
\end{array} \right),$$
has rank $n-1$. As in Lemma \ref{GenericProj}, we can assume that the minor 
$$\left\vert \begin{array}{ccc}
\frac{\partial f_1}{\partial x_1} (x) & \cdots & \frac{\partial f_1}{\partial x_{n-1}} (x)\cr
\vdots & \ddots & \vdots \cr
\frac{\partial f_{c_S}}{\partial x_1}(x) & \cdots & \frac{\partial f_{c_S}}{\partial x_{n-1}} (x) \cr
\frac{\partial m_{c_S+k} }{\partial x_1} (x) & \cdots & \frac{\partial m_{c_S+k}}{\partial x_{n-1}} (x)\cr
\vdots & \ddots & \vdots \cr
\frac{\partial m_n}{\partial x_1}(x) & \cdots & \frac{\partial m_n}{\partial x_{n-1}} (x) \cr
v_1^1 & \cdots & v_1^{n-1} \cr
\vdots & \ddots & \vdots \cr
v_{k-2}^1 & \cdots & v_{k-2}^{n-1} \cr  
\end{array} \right\vert,$$
does not vanish. The set  $N_1$ is defined locally by the vanishing of the $f_i$'s, the $m_j$'s and the determinant of $C$. Because the above $(n-1) \times (n-1)$-minor is not zero, we see that the gradient vectors of the functions defining $N_1$ have rank $n-k+2$ and $N_1$ is a submanifold of dimension $nk+k-2$. In the same way, the sets $N_i$, $i \ge 2$, have dimension less than $nk+k-2$ and $N$ is a definable set of dimension at most $nk+k-2$. Let us consider the projection $\rho : N \rightarrow \mathcal{V} \setminus \overline{\bf D}, (x,V ) \mapsto V$. Let ${\bf D}'$ be the union of the discriminants of the restrictions of $\rho$ to the $N_i$'s. It is a definable set of dimension strictly less than $nk$. Hence for all $V \in \mathcal{V} \setminus (\overline{\bf D} \cup {\bf D}')$, $\rho^{-1}(V)$ is a definable set of dimension $k-2$ (possibly empty). But $\rho^{-1}(V)$ is exactly $\Sigma_P^S \setminus (\Sigma_P^{S,0} \cup \Sigma_P^{'S})$, where $V=(v_1,\ldots,v_k)$ and $P=[v_1,\ldots,v_k ]$.  We take $Z_P^S = \Sigma_P^{'S} \cup \rho^{-1}(V)$. It is closed because $\Sigma_P^{'S}$ is closed and $\overline{\rho^{-1}(V)} \subset \rho^{-1}(V) \cup \Sigma_P^{'S}$. 

For $k=d_S+1$, the statement is just a reformulation of Lemma \ref{GenericProj}. \endproof

\begin{lemma}\label{GenericPolarImage}
Let $S$ be a stratum of $X$ such that $d_S< n$. If $k \le d_S+1$ then for almost all $P \in G_n^k$, $\Delta_P^S$ is a definable set of dimension $k-1$ (or empty).
\end{lemma}
\proof  The set $\Delta_P^S$ is definable  as the projection of a definable set. Moreover ${\rm dim}(\Delta_P^S) \le {\rm dim}(\Sigma_P^S)=k-1$. The projection ${\pi_P^S}_{\vert \Sigma_P^{S,0}} : \Sigma_P^{S,0} \rightarrow \Delta_P^S$ is a local diffeomorphism because 
for all $x \in \Sigma_P^{S,0}$, $T_x  \Sigma_P^S \cap P^\perp = \{0\}$. Let $x$ be a point in $\Sigma_P^{S,0}$. Then  there exists a neighborhood $U_x$ of $x$ included in $\Sigma_P^{S,0}$ on which ${\pi_P^S}_{\vert \Sigma_P^{S,0}} $ is a diffeomorphism and so ${\rm dim} \ \pi_P^S (U_x)=k-1$. But $\pi_P^S(U_x)$ is included in $\Delta_P^S$. Therefore the dimension of  $\Delta_P^S$ is greater or equal to $k-1$. \endproof
Let $D({\pi_P^S}_{\vert \Sigma_P^{S,0}})$ be the set of double-points of ${\pi_P^S}_{\vert \Sigma_P^{S,0}}$, i.e.,
$$D({\pi_P^S}_{\vert \Sigma_P^{S,0}}) = \left\{ x \in \Sigma_P^{S,0} \ \vert \ \exists y \in \Sigma_P^{S,0} \hbox{ with } x \not=  y \hbox{ and } \pi_P(x)=\pi_P(y) \right\}.$$

\begin{lemma}\label{GenericDoublePoints}
Let $S$ be a stratum of $X$ such that $d_S< n$. If $k \le d_S+1$ then for almost all $P \in G_n^k$, $D({\pi_P^S}_{\vert \Sigma_P^{S,0}})$ is a definable set of dimension at most $k-2$ (or empty).
\end{lemma}
\proof The set $D({\pi_P^S}_{\vert \Sigma_P^{S,0}})$ is clearly definable. Let us treat first the case $k \le d_S$. We keep the notations of the previous lemmas. Let $L$ be the following definable set:
$$\displaylines{
L = \Big\{ (x,y,v_1,\ldots,v_k) \in \mathbb{R}^n \times \mathbb{R}^n \times (\mathcal{V} \setminus \overline{{\bf D} \cup {\bf D'}}) \ \vert \ x \in \Sigma_P^{S,0}, y \in \Sigma_P^{S,0}, x \not=y  \hfill \cr
\hfill  \hbox{and } \langle v_j,x \rangle= \langle v_j, y \rangle \hbox{ for all } j \in \{1,\ldots,k \} \hbox{ where } P= [ v_1,\ldots, v_k ]  \Big\}. \quad \cr
}$$
This set is a definable submanifold of dimension $nk+k-2$. Let us explain this assertion. As explained in Lemma \ref{GenericProj}, if $x \in \Sigma_P^{S,0}$ then there exist smooth definable functions $f_1,\ldots,f_{c_S}$ and $m_{c_S+k}(-,V),\ldots,m_n(-,V)$, where $V=(v_1,\ldots,v_k)$, such that in a neighborhood of $x$, $\Sigma_P^{S,0}$ is the set $$\{f_1=\cdots=f_{c_S}=m_{c_S+k}(-,V)=\cdots=m_n(-,V)=0\},$$ and the gradient vectors of these functions are linearly independent. Furthermore, since $x$ is a fold point, 
$${\rm rank}(\nabla f_1(x), \ldots, \nabla f_{c_S}(x), \nabla_x m_{c_S+k}(x,V), \ldots, \nabla_x m_n(x,V), v_1,\ldots,v_k)=n.$$
Similarly, there exist smooth definable functions $g_1,\ldots,g_{c_S}$ and $\tilde{m}_{c_S+k}(-,V),\ldots,$ $\tilde{m}_n(-,V)$ such that in a neighborhood of $y$, $\Sigma_P^{S,0}$ is the set 
$$\{g_1=\cdots=g_{c_S}=\tilde{m}_{c_S+k}(-,V)=\cdots=\tilde{m}_n(-,V)=0\},$$ and the gradient vectors of these functions are linearly independent. We also have that
$${\rm rank}(\nabla g_1(y), \ldots, \nabla g_{c_S}(y), \nabla_y \tilde{m}_{c_S+k}(y,V), \ldots, \nabla_y \tilde{m}_n(y,V), v_1,\ldots,v_k)=n.$$
Therefore $L$ is locally given by the vanishing of the functions 
$$f_1,\ldots,f_{c_S},m_{c_S+k},\ldots,m_n,g_1,\ldots,g_{c_S}, \tilde{m}_{c_S+k},\ldots,\tilde{m}_n,$$ and $$\langle v_1,x-y\rangle, \ldots, \langle v_k, x-y \rangle.$$ Let us remark that the $f_i$'s depend on $x$, the $m_j$'s on $x$ and $V$, the $g_i$'s on $y$ and the $\tilde{m}_j$'s on $y$ and $V$. We will show that the gradient vectors of these functions have rank $2n-k+2$. As in the proof of Lemma \ref{GenericProj}, we can assume that $v_k$ belongs to $[ \nabla f_1 (x),\ldots,\nabla f_{c_S}(x), v_1,\ldots,v_{k-1} ]$. Since $x \not= y$, there is $j \in \{1,\ldots,n\}$ such that $x_j-y_j \not= 0$. Let us consider the following matrix $J$ of size $(2n-k+2) \times (2n+1)$:
$$J = \tiny{ \left( \begin{array}{ccc}
 \nabla f_1 (x)  & {\bf 0}  & 0\cr
 \vdots &   \vdots & \vdots  \cr
 \nabla f_{c_S} (x)  &  {\bf 0} & 0\cr
\nabla_x m_{c_S+k} (x,V) & {\bf 0}  &  * \cr
\vdots & \vdots &  \vdots  \cr
\nabla_x m_{n} (x,V) & {\bf 0} &   * \cr
 {\bf 0} &  \nabla g_1 (y) & 0\cr
 \vdots & \vdots & \vdots  \cr
 {\bf 0}&  \nabla g_{c_S} (y) &  0\cr
 {\bf 0} & \nabla_y \tilde{m}_{c_S} (y,V) &  *\cr
 \vdots & \vdots &  \vdots  \cr 
 {\bf 0} & \nabla_y \tilde{m}_{n} (y,V)  & *\cr
 v_1   & -v_1 &  0\cr
 \vdots & \vdots &  \vdots  \cr 
 v_{k-1} & -v_{k-1} &  0\cr
 v_{k}  & -v_{k} & x_j-y_j \cr
\end{array} \right),}$$
where ${\bf 0}=(0,\ldots,0)$ in $\mathbb{R}^n$.
It is a submatrix of the matrix whose lines are the gradient vectors of the $f_i$'s, the $m_j$'s, the $g_i$'s, the $\tilde{m}_j$'s and the $k$ functions $\langle v_i, x-y \rangle$.  Since $v_k$ belongs to $[ \nabla f_1 (x),\ldots,\nabla f_{n-d_S}(x), v_1,\ldots,v_{k-1} ]$, $J$ has the same rank as the matrix $\tilde{J}$ given by
$$\tilde{J} =  {\tiny \left( \begin{array}{ccc}
 \nabla f_1 (x)  & {\bf 0}  & 0\cr
 \vdots &   \vdots & \vdots  \cr
 \nabla f_{c_S} (x)  &  {\bf 0} & 0\cr
\nabla_x m_{c_S+k} (x,V) & {\bf 0}  &  * \cr
\vdots & \vdots &  \vdots  \cr
\nabla_x m_{n} (x,V) & {\bf 0} &   * \cr
 {\bf 0} &  \nabla g_1 (y) & 0\cr
 \vdots & \vdots & \vdots  \cr
 {\bf 0}&  \nabla g_{c_S} (y) &  0\cr
 {\bf 0} & \nabla_y \tilde{m}_{c_S} (y,V) &  *\cr
 \vdots & \vdots &  \vdots  \cr 
 {\bf 0} & \nabla_y \tilde{m}_{n} (y,V)  & *\cr
 v_1   & -v_1 &  0\cr
 \vdots & \vdots &  \vdots  \cr 
 v_{k-1} & -v_{k-1} &  0\cr
{\bf 0}  & w & x_j-y_j \cr
\end{array} \right),}$$
where $w=(w_1,\ldots,w_n)$ is a linear combination of the $-v_i$'s and thus belongs to $P$. Let us denote the lines of $\tilde{J}$ by 
$${\bf F}_1,\ldots,{\bf F}_{c_S}, {\bf M}_{c_S+k},\ldots,{\bf M}_n,{\bf G}_1,\ldots,{\bf G}_{c_S}, {\bf \tilde{M}}_{c_S+k},\ldots,{\bf \tilde{M}}_n,$$  and $${\bf V}_1,\ldots,{\bf V}_k,$$ and show that they are linearly independent. 
So let us suppose that
$$\sum_i \alpha_i {\bf F}_i + \sum_j \beta_j {\bf M}_j + \sum_i \gamma_i {\bf G}_i + \sum_j \delta_j {\bf \tilde{M}}_j + \sum_i \xi_i {\bf V}_i=0.$$
Since $${\rm rank}(\nabla f_1(x),\ldots,\nabla f_{c_S}(x),\nabla_x m_{c_S+k}(x,V),\ldots, \nabla_x m_n (x,V), v_1,\ldots,v_{k-1})=n,$$ we see that $$\alpha_1=\cdots=\alpha_{c_S}=\beta_{c_S+k}=\cdots=\beta_n=\xi_1=\cdots=\xi_{k-1}=0.$$ We are left with the equality
$$ \sum_i \gamma_i {\bf G}_i + \sum_j \delta_j {\bf \tilde{M}}_j + \xi_k {\bf V}_k=0,$$
which implies that $\xi_k w$ belongs to 
$$[\nabla g_1 (y),\ldots,\nabla g_{c_S}(y),\nabla_y \tilde{m}_{c_S+k}(y,V),\ldots, \nabla_y \tilde{m}_n (y,V)].$$ But $y$ is a fold point so ${\rm dim}(T_y S \cap P^\perp)=d_S-k+1$ and ${\rm dim}(T_y \Sigma_P^S \cap P^\perp)=0$. Therefore, ${\rm dim}(N_y S \cap P)={\rm dim}(N_y \Sigma_P^S \cap P) =1$ and $N_y S \cap P= N_y \Sigma_P^S \cap P$. So the vector $\xi_k w$ belongs to $[ \nabla g_1(y),\ldots,\nabla g_{c_S}(y) ]$, and $\delta_{c_S+k}=\cdots=\delta_n=0$. Since $x_j-y_j \not= 0$, we find that $\xi_k=0$ and finally that $\gamma_1=\cdots=\gamma_{c_S}=0$. We conclude that the gradient vectors of the functions which define $L$ locally at $(x,y,V)$ are linearly independent and that $L$ has dimension $nk+k-2$. As in the previous statements, we see that for almost all $(v_1,\ldots,v_k)$ in $\mathcal{V} \setminus \overline{{\bf D} \cup {\bf D'}}$, the set
$$\displaylines{
\quad \Big\{ (x,y) \in \mathbb{R}^n \times \mathbb{R}^n \ \vert \ x \in \Sigma_P^{S,0}, y \in \Sigma_P^{S,0}, x \not=y  \hfill \cr
\hfill \hbox{and } \pi_P^S(x)=\pi_P^S(y) \hbox{ where } P=[ v_1,\ldots,v_k ] \Big\}, \quad \cr
}$$
has dimension at most $k-2$. Therefore $D({\pi_P^S}_{\vert \Sigma_P^{S,0}})$ has dimension at most $k-2$ because it is the image of this set by the projection $(x,y) \mapsto x$. 

Let us treat the case $k=d_S+1$. Here we recall that $\Sigma_P^{S,0}$ is the set of regular point of $\pi_P^S$. Let $Q$ be the following definable set:
$$\displaylines{
Q = \Big\{ (x,y,v_1,\ldots,v_k) \in \mathbb{R}^n \times \mathbb{R}^n \times \mathcal{V}  \ \vert \ x \in S, y \in S, x \not=y  \hfill \cr
\hfill  \hbox{and } \langle v_j,x \rangle= \langle v_j, y \rangle \hbox{ for all } j \in \{1,\ldots,k \} \hbox{ where } P= [ v_1,\ldots, v_k ]  \Big\}. \quad \cr
}$$
This set is a definable submanifold of dimension $nk+k-2$. The proof of this fact is an in the previous case. Locally at a point $(x,y,V)$, $Q$ is given by the vanishing of functions $f_1,\ldots,f_{c_S},g_1,\ldots,g_{c_S}$ and the functions $\langle v_1, x-y \rangle,\ldots, \langle v_k, x-y \rangle$. As before the $f_i$'s depend only on $x$ and the $g_i$'s only on $y$. It is not difficult to see that the gradient vectors of these functions are linearly independent because $x-y \not= 0$. Therefore for almost all $V$ in $\mathcal{V}$, the set 
$$\Big\{ (x,y) \in S \times S \ \vert \ x \not= y \hbox{ and } \pi_P^S(x)=\pi_P^S(y) \hbox{ where } P= [ v_1,\ldots, v_k ]  \Big\},$$
has dimension at most $k-2$ and so has $D({\pi_P^S}_{\vert \Sigma_P^{S,0}})$.
\endproof

\begin{remark}
{\rm As noticed to the author by Terry Gaffney, we can also use Mather's technology of modular submanifolds to prove the above lemma.
Double points of fold singularities form a contact class and a Thom-Boardmann manifold. Hence by \cite{Mather}, pages 233 and 234, they form a modular submanifold. Therefore by Theorem 1 in \cite{Mather}, for almost all $P \in G_n^k$, $\pi_P^S$ is transverse with respect to this submodular manifold. 
Our proof avoids the notions of modular manifolds and contact classes and enables us to stay in the definable setting.}
\end{remark}

For each stratum $S$ and each $P \in G_n^k$, $ k \le d_S+1$, let $\Sigma_{P,{\rm lim}}^S$ be the following set:
$$\displaylines{
\quad \Sigma_{P,{\rm lim}}^S = \left\{ 
x \in \bar{S} \setminus S \ \vert \ \exists (x_m)_{m \in \mathbb{N}}  \hbox{ in } S \hbox{ such that } x_m \rightarrow x , T_{x_m} S \rightarrow T  \right. \hfill \cr
   \hfill \left. \hbox{and dim}(T\cap P^\perp) \ge d_S-k+1 \right\}. \quad \cr
}$$

\begin{proposition}\label{GenericProjLim}
Let $S$ be  a stratum of $X$. If $k \le d_S+1$ then
for almost all $P \in G_n^k$, $\Sigma^S_{P,{\rm lim}}$ is a definable set of dimension at most $k-2$ (or empty).
\end{proposition}
\proof We treat first the case $d_S <n$ and $k \le d_S$. In Lemma \ref{GenericProj}, we proved that  the following definable set $M$:
$$\displaylines{
M = \Big\{ (x,v_1,\ldots,v_k) \in \mathbb{R}^n \times \mathcal{V} \ \vert \  x \in S \hbox{ and dim }(T_x S \cap P^\perp) \ge d_s-k+1  \hfill \cr
 \hfill    \hbox{ where } P=[ v_1,\ldots, v_k ] \Big\}.  \qquad \cr
 }$$
had  dimension $nk+k-1$. 
 
 If $x \in \Sigma^S_{P,{\rm lim}}$ then there exists a sequence $(x_m)_{m \in \mathbb{N}} $ in $S$ such that $x_m \rightarrow x$, $T_{x_m} \rightarrow T$ and ${\rm dim}(T \cap P^\perp) \ge d_S -k +1$. This implies that there exists a sequence $(P_m)_{m \in \mathbb{N}} $ in $G_n^k$ such that $P_m \rightarrow P$ and ${\rm dim}(T_{x_m} S \cap P_m^\perp) \ge d_S -k +1$, because we can find a subspace $L_m$ in $T_{x_m} S$ such that $L_m \rightarrow T \cap P^\perp$. Conversely if $T_{x_m} S \rightarrow T$, $P_m \rightarrow P$ and ${\rm dim}(T_{x_m} S\cap P_m^\perp) \ge d_S -k +1$ then ${\rm dim}(T \cap P^\perp) \ge d_S -k +1$. Therefore $x \in \Sigma_{P,{\rm lim}}^S$ if and only if $(x,v_1,\ldots,v_k) \in \overline{M} \setminus M$ where $P=[ v_1,\ldots,v_k ]$. 

The set $\overline{M} \setminus M$ is definable of dimension $nk+k-2$. We can conclude that for almost all $V=(v_1,\ldots,v_k)$ in $\mathcal{V}$, $\pi^{-1}(V) \cap \overline{M} \setminus M$ has dimension at most $k-2$, where $\pi$ is the projection $(x,V)  \mapsto V$.

If $d_S <n$ and $k=d_S+1$ then $\Sigma_{P,\lim}^S$ is just $\overline{S} \setminus S$, which has dimension at most $d_S-1=k-2$. 

If $d_S=n$, the result is obvious because $\Sigma^S_{P,{\rm lim}}$ is empty. \endproof

Let us fix now a stratum $S$ such that $d_S <n$ and consider the strata that contain $S$ in their frontier. We denote them by $S_1,\ldots,S_r$. Note that for $i=1,\ldots,r$, $\Sigma_{P,\lim}^{S_i} \cap S \subset \Sigma_P^S $ by Whitney condition (a) (see the proof of Lemma \ref{Transverse0}). Applying the previous results, we see that for almost all $P \in G_n^k$, the  set $T_P^S$ defined by 
$$ T_P^S =Z_P^S \bigcup \overline{D({\pi_P^S}_{\vert \Sigma_P^{S,0}})}  \bigcup \overline{( \cup_{i=1}^r \Sigma_{P,\lim}^{S_i}) \cap S} ,$$ is a closed definable set of dimension at most $k-2$.  

From now on, we fix such a generic $P$ in $G_n^k$, $k=1,\ldots,d_S+1$. Let $y$ be a point in $\Delta_P^S \setminus (\pi_P^S (T_P^S) \cup {\rm Sing}(\Delta_P^S))$. Since $y$ does not belong to $\pi(\overline{D({\pi_P^S}_{\vert \Sigma_P^{S,0}})})$, there is a unique $x$ in 
$$\Sigma_P^{S,0} \setminus \left[ \overline{D({\pi_P^S}_{\vert \Sigma_P^{S,0}})}  \bigcup \overline{( \cup_{i=1}^r \Sigma_{P,\lim}^{S_i})\cap S} \right],$$
such that $\pi_P^S(x)=y$. Let $u$ be a unit vector in $P$ not belonging to $T_y \Delta_P^S$. Let $Q_u= P^\perp \oplus u$ and let $Q_{u,x}$ be the affine space through $x$ parallel to $Q_u$. 
Since $T_x \Sigma_P^S \oplus P^\perp$ has dimension $n-1$ and $u$ does not belong to this vector space, $Q_{u,x}$ intersects $\Sigma_P^S$ and $S$ transversally at $x$.

\begin{lemma}\label{MorseIndex1}
For any unit vector $u$ in $P \setminus T_y \Delta_P^S$, the point $x$ is a non-degenerate critical point of $u^*_{\vert Q_{u,x} \cap S}$. Furthermore if $\nu$ is the unit normal vector to $\Delta_P^S$ at $y$ in $P$ such that $\langle u, \nu \rangle >0$, then $u^*_{\vert Q_{u,x} \cap S}$ and $\nu^*_{\vert Q_{\nu,x} \cap S}$ have the same tangential Morse index.
\end{lemma}
\proof Let us treat first the case $k \le d_S$. It is clear that $x$ is a critical point of $u^*_{\vert Q_{x,u} \cap S}$ because $(P^\perp \oplus u ) \cap T_x S = P^\perp \cap T_x S$.  
We keep the notations of the previous proofs. For simplicity, we assume that $y=0$ in $P$. 

Note that $\nu$ is a normal vector to $S$ at $x$ because $T_x S = T_x \Sigma_P^S \oplus (T_x S \cap P^\perp)$.  Let $(v_1,\ldots,v_{k-1})$ be an orthonormal basis of $T_y \Delta_P^S$ such that $(\nu,v_1,\ldots,v_{k-1})$ is a positive orthonormal basis of $P$. Then $Q_{\nu,x} \cap S$ is the manifold $$\{f_1=\cdots=f_{c_S}=v_1^*=\cdots=v_{k-1}^*=0 \},$$ in a neighborhood of $x$. We can assume that the following minor:
$$m^\nu =\frac{\partial(v_1^*,\ldots,v_{k-1}^*,f_1,\ldots,f_{c_S})}{\partial(x_1,\ldots,x_{k-1},x_k, \ldots, x_{c_S+k-1})},$$
does not vanish at $x$. Therefore $x$ is a critical point of $\nu^*_{\vert Q_{\nu,x} \cap S}$ if and only if the following minors:
$$m^\nu_j =\frac{\partial(\nu^*,v_1^*,\ldots,v_{k-1}^*,f_1,\ldots,f_{c_S})}{\partial(x_1,\ldots,x_{k-1},x_k, \ldots, x_{c_S+k-1},x_j)},$$
$j=c_S+k,\ldots,n$, vanish at $x$. But as explained in the previous lemmas, $\Sigma_P^S$ is the set  $\{f_1=\cdots=f_{c_S}=m^\nu_{c_S+k}=\cdots=m^\nu_{n}=0 \}$ in a neighborhood of $x$. Since $\{v^*_1=\cdots=v^*_{k-1}=0 \}$ intersects $\Sigma_P^S$ transversally at $x$, the following $n \times n$ minor:
$$ \frac{\partial(v_1^*,\ldots,v_{k-1}^*,f_1,\ldots,f_{c_S},m^\nu_{c_S+k},\cdots,m^\nu_{n})}{\partial(x_1,\ldots,x_n)},$$
does not vanish at $x$. By the computations of Szafraniec (\cite{Szafraniec}, pages 248-250), this exactly means that $\nu^*_{\vert Q_{\nu,x} \cap S}$ has a Morse critical point. Moreover the tangential index of $\nu^*_{\vert Q_{\nu,x} \cap S}$ at $x$ is the sign of 
$$\displaylines{
\qquad (-1)^{(c_S+k-1)(d_S-k+1)} (m^\nu (x) )^{d_S+k+2} \times \hfill \cr
\hfill  \frac{\partial(v_1^*,\ldots,v_{k-1}^*,f_1,\ldots,f_{c_S},m^\nu_{c_S+k},\cdots,m^\nu_{n})}{\partial(x_1,\ldots,x_n)}(x). \qquad \cr
}$$
Let $(w_1,\ldots,w_{k-1})$ be an orthonormal basis of $u^*$ such that $(u,w_1,\ldots,w_{k-1})$ is a positive orthonormal basis of $P$. There exists a positive orthogonal $k \times k$ matrix $A=(a_{ij})$ such that 
$$\left\{ \begin{array}{c}
u = a_{11} \nu + \cdots + a_{1k} v_{k-1} \cr
\vdots \cr
w_{k-1}= a_{k1} \nu + \cdots + a_{kk} v_{k-1}, \cr
\end{array} \right.$$
where $a_{11} >0$.  Then $Q_{u,x} \cap S$ is the manifold $$\{f_1=\cdots=f_{c_S}=w_1^*=\cdots=w_{k-1}^*=0 \},$$ in a neighborhood of $x$. Let $m^u$ be the following minor:
$$m^u =\frac{\partial(w_1^*,\ldots,w_{k-1}^*,f_1,\ldots,f_{c_S})}{\partial(x_1,\ldots,x_{k-1},x_k, \ldots, x_{c_S+k-1})}.$$
Since $\nu$ belongs to $[\nabla f_1 (x), \ldots, \nabla f_{c_S}(x) ]$ and $A$ is orthogonal, it is not difficult to see that $m^u(x)= a_{11} m^{\nu}(x)$ and thus does not vanish. For $j=c_S+k,\ldots,n$, let $m_j^u$ be the following minor:
$$m^u_j =\frac{\partial(u^*,w_1^*,\ldots,w_{k-1}^*,f_1,\ldots,f_{c_S})}{\partial(x_1,\ldots,x_{k-1},x_k, \ldots, x_{c_S+k-1},x_j)}.$$
Since $A$ is orthogonal, $m_j^u$ is equal to $m_j^{\nu}$ and since $\nu$ belongs to $[\nabla f_1 (x), \ldots, \nabla f_{c_S}(x) ]$,
$$\displaylines{
\quad \frac{\partial(w_1^*,\ldots,w_{k-1}^*,f_1,\ldots,f_{c_S},m^u_{c_S+k},\cdots,m^u_{n})}{\partial(x_1,\ldots,x_n)}(x)  =  \hfill \cr
\hfill a_{11} \frac{\partial(v_1^*,\ldots,v_{k-1}^*,f_1,\ldots,f_{c_S},m^\nu_{c_S+k},\cdots,m^\nu_{n})}{\partial(x_1,\ldots,x_n)}(x). \quad \cr
}$$
We can conclude that $u^*_{\vert Q_{u,x} \cap S}$ has a Morse critical point at $x$ and that the two functions have the same tangential Morse index because $a_{11} >0$.

The case $k=d_S+1$ is easy because in this situation, $$Q_{u,x} \cap S = Q_{\nu,x} \cap S = \{x\},$$ in a neighborhood of $x$ and the tangential Morse index is $1$.
\endproof

By Whitney condition (a), $Q_{u,x}$ intersects also the strata that contains $x$ in their frontier transversally in a neighborhood of $x$ and hence, $Q_{u,x}\cap X$ is Whitney stratified around $x$.
\begin{lemma}
For any unit vector $u$ in $P \setminus T_y \Delta_P^S$, the function $u^*_{\vert Q_{u,x} \cap X}$ is a stratified Morse function in a neighborhood of $x$, with a critical point at $x$.
\end{lemma}
\proof We already know that $u^*_{\vert Q_{u,x} \cap S}$ has a Morse critical point at $x$. The fact that $x$ does not belong to $\cup_{i=1}^r \Sigma_{P,\lim}^{S_i}$ implies that $P^\perp$ intersects the strata $S_1,\ldots,S_r$ transversally in a neighborhood of $x$. Since $Q_{u,x}$ is also transverse to these strata around $x$, $u^*_{\vert Q_{u,x} \cap X}$ has no critical point on these strata. 

Let us check now that $u$ is not perpendicular to any limit tangent space at $x$ in $Q_{u,x} \cap X$. If it is not the case then there is a sequence of points $(x_m)_{m \in \mathbb{N}}$ in a stratum $S_i$ such that $x_m \rightarrow x$, and $T_{x_m} S_i$ tends to $T$ where $T$ satisfies the following condition: $u \perp (T \cap Q_{u,x})$ in $Q_{u,x}$. This implies that $T \cap P^\perp = T \cap Q_{u,x}$ and that ${\rm dim}(T \cap P^\perp) \ge d_{S_i}-k +1$, which contradicts the fact that $x$ is not in $\Sigma_{P,\lim}^{S_i}$.   \endproof

\begin{lemma}\label{MorseIndex2}
Let $u$ be a  unit vector in $P \setminus T_y \Delta_P^S$ and let $\nu$ be the unit normal vector to $\Delta_P^S$ at $y$ in $P$ such that $\langle u, \nu \rangle >0$, then $u^*_{\vert Q_{u,x} \cap X}$ and $\nu^*_{\vert Q_{\nu,x} \cap X}$ have the same stratified Morse index.
\end{lemma}
\proof Let us study the behavior of $\pi_P^X : X \rightarrow P$ in the neighborhood of $x$. Since $P^\perp$ is transverse to the strata $S_1,\ldots,S_r$ in a neighborhood of $x$ and since ${\pi_P^S}_{\vert \Sigma_P^{S,0}}$ is a local diffeomorphism, we see that $\pi_P^{-1}(y)$ intersects $X \setminus \{x\}$ transversally in a neighborhood of $x$. By the Curve Selection Lemma applied to the stratified space $\pi_P^{-1}(y) \cap X$, the function $\omega_x$ is a submersion on $\pi_P^{-1}(y)  \cap X \setminus \{x\}$, where $\omega_x(z)= \Vert z-x \Vert^2$. Hence, for $\epsilon >0$ small enough so that $S_\epsilon(x)$ is transverse to $X$ in the stratified sense, there exists a small neighborhood $U_\epsilon$ of $y$ in $P$ such that for all $y'$ in $U_\epsilon$, $\pi_P^{-1}(y')$ is transverse to $X \cap S_\epsilon(x)$ by Lemma \ref{Transverse}. 

Now in a neighborhood of $y$, the discriminant of $\pi_P^X $ is exactly $\Delta_P^S$, so by the Thom-Mather  lemma and shrinking $U_\epsilon$ if necessary, we can say that if $y_1$ and $y_2$ lie in the same component of $U_\epsilon \setminus \Delta_P^S$ then $\pi_P^{-1}(y_1) \cap X \cap B_\epsilon(x)$ and $\pi_P^{-1}(y_2) \cap X \cap B_\epsilon(x)$ are homeomorphic and $\chi(\pi_P^{-1}(y_1) \cap X \cap B_\epsilon(x)) = \chi(\pi_P^{-1}(y_2) \cap X \cap B_\epsilon(x)) $.  The Morse index of $u^*_{\vert Q_{u,x} \cap X}$ at $x$ is 
$$1-\chi (Q_{u,x} \cap X \cap \{u^*= u(x)-\delta' \} \cap B_{\epsilon'}(x)),$$
where $0 < \delta' \ll \epsilon' \ll 1$, that is 
$$1- \chi( \pi_P^{-1}(y-\delta' u ) \cap X \cap B_{\epsilon'}(x)).$$
Similarly the Morse index of $\nu^*_{\vert Q_{\nu,x} \cap X}$ at $x$ is
$$1-\chi (\pi_P^{-1}(y-\delta" \nu) \cap X \cap B_{\epsilon"}(x)),$$
where $0 < \delta" \ll \epsilon" \ll 1$. We can choose $\epsilon$ small enough so that the Morse index of $u^*_{\vert Q_{u,x} \cap X}$ at $x$ is 
$$1- \chi( \pi_P^{-1}(y-\delta' u ) \cap X \cap B_\epsilon(x)).$$
where $0 < \delta' \ll \epsilon$  and $y-\delta' u \in U_\epsilon \setminus \Delta_P^S$, and  the Morse index of $\nu^*_{\vert Q_{\nu,x} \cap X}$ at $x$ is
$$1-\chi (\pi_P^{-1}(y-\delta" \nu) \cap X \cap B_{\epsilon}(x)),$$
where $0 < \delta" \ll \epsilon $ and $y-\delta" \nu \in U_\epsilon \setminus \Delta_P^S$. But, if $\delta'$ and $\delta"$ are small enough, then $y-\delta' u$ and $y-\delta" \nu$ lie in the same component of $U_\epsilon \setminus \Delta_P^S$, since $\langle u, \nu \rangle >0$. \endproof

\begin{corollary}
Let $u$ be a unit vector in $P \setminus T_y \Delta_P^S$ and let $\nu$ be the unit normal vector to $T_y \Delta_P^S$ in $P$ such that $\langle u, \nu \rangle >0$, then $u^*_{\vert Q_{u,x} \cap X}$ and $\nu^*_{\vert Q_{\nu,x} \cap X}$ have the same normal Morse index at $x$.
\end{corollary}
\proof Combine Lemma \ref{MorseIndex1}, Lemma  \ref{MorseIndex2} and the results of Subsection 2.2. \endproof

\begin{definition}
{\rm For $P$ generic in the sense of Lemmas \ref{GenericProj}, \ref{GenericFolds}, \ref{GenericDoublePoints} and Proposition \ref{GenericProjLim} and for $y \in \Delta_P^S \setminus ( \pi_P^S(T_P^S) \cap {\rm Sing} (\Delta_P^S))$, we define 
$$\alpha(\pi_P^S,y) = \frac{1}{2} \Big( {\rm ind}(\nu^*, X \cap Q_{\nu,x} ,x) + {\rm ind}(-\nu^*, X \cap Q_{\nu,x} ,x) \Big),$$
where $\nu$ is a unit normal vector to $T_y \Delta_P^S$ and $y=\pi_P^S(x)$.}
\end{definition} 
Note that by Lemma \ref{MorseIndex2}, if $u$ is a unit vector in $P \setminus T_y \Delta_P^S$, then 
$$\alpha(\pi_P^S,y) = \frac{1}{2} \Big( {\rm ind}(u^*, X \cap Q_{u,x} ,x) + {\rm ind}(-u^*, X \cap Q_{u,x} ,x) \Big).$$

\begin{lemma}\label{IndiceLocConstant}
There exists a closed definable set $T'^S_P \subset \Sigma_P^S$ of dimension less or equal to $k-2$ such that for $x$ in $\Sigma_P^S \setminus (T_P^S \cup T'^S_P)$,  there exists a neighborhood $U_x$ of $x$ in $\Sigma_P^S \setminus (T_P^S \cup T'^S_P)$ such that $\alpha(\pi_P^S,\pi_P^S(x))= \alpha(\pi_P^S,\pi_P^S(x'))$ for all $x' \in U_x$.
\end{lemma}
\proof Let $d : \mathbb{R}^n \rightarrow \mathbb{R}$ be the distance function to $\Sigma_P^{S}$. It is a continuous definable function and moreover there exists an open definable neighborhood $\mathcal{U}$ of $\Sigma_P^S \setminus T_P^S$ such that $d$ is smooth on $\mathcal{U} \setminus (\Sigma_P^S \setminus T_P^S)$. By the Curve Selection Lemma, we can also assume that $d$ is a (stratified) submersion on $X \cap [ \mathcal{U} \setminus (\Sigma_P^S \setminus T_P^S)]$.  

Let $x$ be a point in $\Sigma_P^S \setminus T_P^S$ and let $\nu$ be a unit normal vector to $\Delta_P^S$ at $y$ in $P$ where $\pi_P^S(x)=y$. Since $x$ is an isolated point in $Q_{\nu,x} \cap \Sigma_P^S$, 
the function $d_{\vert Q_{\nu,x} }: Q_{\nu,x} \rightarrow \mathbb{R}$ is a distance function to $x$ and therefore
$$ {\rm ind}(\nu^*, X \cap Q_{\nu,x} ,x) =1-\chi \left(\pi_P^{-1} (y-\delta \nu) \cap X \cap d^{-1}([0,\epsilon]) \right),$$
where $0 < \delta \ll \epsilon \ll 1$. 
As explained in Lemma \ref{MorseIndex2}, there exists a neighborhood $W$ of $y$ in $P$ such that for all $y_1 \in W\setminus \Delta_P^S$, $\pi_P^{-1}(y_1)$ intersects $X$ transversally in a neighborhood of $x$. 
Moreover there exists $\epsilon_x >0$ and a neighborhood $\Omega$ of $x$ in $\Sigma_P^S \setminus T_P^S$ such that for all $z$ in $\Omega$, $\pi_P^{-1}( \pi_P(z)) $ intersects $X \setminus \{z\}$ transversally in $\{ 0 < d \le \epsilon_x\}$. If it is not the case, we can find a sequence of points $(z_m)_{m \in \mathbb{N}}$ in $\Sigma_P^S \setminus T_P^S$ that tends to $x$ and a sequence of non-negative real numbers $(\epsilon_m)_{m \in \mathbb{N}}$ that tends to $0$ such that $\pi_P^{-1} (\pi_P(z_m))$ does not intersect $X \setminus \{ z_m \}$ transversally in $\{ 0 < d \le \epsilon_m\}$. This implies that $x$ belongs to $\cup_{i=1}^r \Sigma_{P,{\rm lim}}^{S_i}$, which is excluded.

As in Lemma \ref{MorseIndex2}, there exists $\epsilon'_x >0$ such that for $0 < \epsilon \le \epsilon'_x$, $\pi_P^{-1}(y)$ intersects $X \cap d^{-1}(\epsilon)$ transversally. We are going to explain that for $x$ generic, we can choose $\epsilon'_x$ uniformly in a neighborhood of $x$. Let us consider the mapping
$$\begin{array}{ccccc}
(d,\pi_P^X) & : & ( X \cap \mathcal{U}) \setminus \Sigma_P^X & \rightarrow & \mathbb{R} \times P \cr
    &   & x  & \mapsto & (d(x),\pi_P^X(x)), \cr
\end{array}$$
where $\Sigma_P^X$ is the union of the $\Sigma_P^{S'}$'s, for $S'$ a stratum of $X$.  Let us call ${\bf D}$ its discriminant. It is a definable set of dimension at most $k$. In $\mathbb{R} \times P$, let $\pi_2$ be the projection $(r,y) \mapsto y$. For $y'$ in $\Delta_P^S \setminus \pi_P^S(T_P^S)$, $\pi_2^{-1}(y') \cap {\bf D}$ is included in the set 
$$\left\{ (r,y') \ \vert \ r \hbox{ critical value of  } d_{\vert \pi_P^{-1}(y') \cap [ (X \cap \mathcal{U}) \setminus \Sigma_P^X]} \right\},$$
and so $\pi_2^{-1}(y') \cap {\bf D}$ has dimension zero. Therefore the set $\pi_2^{-1} (\Delta_P^S \setminus \pi_P^S(T_P^S)) \cap {\bf D}$ has dimension at most $k-1$ and the set
$$K_P^S = \overline{\pi_2^{-1} (\Delta_P^S \setminus \pi_P^S(T_P^S)) \cap {\bf D}} \cap (\Delta_P^S \setminus \pi_P^S(T_P^S)),$$
has dimension at most $k-2$. Let us set $T'^S_P= \overline{(\pi_P^S)^{-1} (K_P^S)}$. Since ${\pi_P^S}_{\vert \Sigma_P^S \setminus T_P^S}$ is a local diffeomorphism, ${\rm dim}(T'^S_P)  \le k-2$.

Now we suppose that $x$ belongs to $\Sigma_P^S \setminus (T_P^S \cup T'^S_P)$. Then there exist a neighborhood $\Omega' \subset \Omega$ of $x$ in $\Sigma_P^S  \setminus (T_P^S \cup T'^S_P)$ and $0< \epsilon'_x \le \epsilon_x$ such that for $0 < \epsilon \le \epsilon'_x$ and $z \in \Omega'$, $\pi_P^{-1}(\pi_P(z))$ intersects $X \cap d^{-1}(\epsilon)$ transversally. If it is not the case, we can find a sequence of points $(z_m)_{m \in \mathbb{N}}$ in $\Omega$ that tends to $x$ and a sequence of non-negative real numbers $(\epsilon_m)_{m \in \mathbb{N}}$ that tends to $0$ such that $\pi_P^{-1} (\pi_P(z_m))$ does not intersect $X \cap d^{-1}(\epsilon_m)$ transversally. But $(\epsilon_m, \pi_p(z_m))$ belongs to $\pi_2^{-1} (\Delta_P^S \setminus \pi_P^S(T_P^S)) \cap {\bf D}$ and so $y=\pi_P^S(x)$ belongs to $K_P^S$, which is excluded.

Let $\Omega"$ be an open neighborhood of $x$ such that $\Omega" \subsetneq \Omega'$. Then for all $0 < \epsilon  \le \epsilon'_x$, there exists $\delta_\epsilon >0$ such that for $z \in \overline{\Omega"}$ and $0 < \delta \le \delta_{\epsilon}$, $\pi_P^{-1}( \pi_P(z)-\delta \nu (z)) $ intersects $X \cap d^{-1}(\epsilon)$ transversally, where $\nu(z)$ is the unit normal vector to $\Delta_P^S$ at $\pi_P(z)$ in $P$ that lies in the same component as $\nu$. If it is not true, then we can find  $\epsilon$ with $0< \epsilon \le \epsilon'$, a sequence of points $(z_m)_{m \in \mathbb{N}}$ in $\overline{\Omega"}$ that tends to $z$ and a sequence of non-negative real numbers $(\delta_m)_{m \in \mathbb{N}}$ that tends to $0$ such that $\pi_P^{-1}( \pi_P(z_m)-\delta_m \nu (z_m)) $ does not intersect $X \cap d^{-1}(\epsilon)$ transversally. But $\pi_P^{-1} (\pi_P(z))$ intersects $X \cap d^{-1}(\epsilon)$ transversally and so for $w$ in a neighborhood of $\pi_P(z)$ in $P$, $\pi_P^{-1}(w)$ intersects $X \cap d^{-1}(\epsilon)$ transversally by Lemma \ref{Transverse}.  

We take $U_x = \Omega" \cap \pi_P^{-1} (W) \cap \Sigma_P^S$ and conclude using the Thom-Mather isotopy lemma as in Lemma  \ref{MorseIndex2}.
 \endproof

\begin{remark}\label{IndiceNormalLocConstant}
{\rm We have proved that the index $ {\rm ind}(\nu^*, X \cap Q_{\nu,x} ,x)$ was locally constant. But since the tangential index $ {\rm ind}_{\rm tg} (\nu^*, X \cap Q_{\nu,x} ,x)$ is locally constant,  the normal index $ {\rm ind}_{\rm nor} (\nu^*, X \cap Q_{\nu,x} ,x)$ is locally constant along $\Sigma_P^S \setminus  (T_P^S \cup T'^S_P)$ as well.}
\end{remark}

\begin{lemma}\label{DecompDiscriminant}
There exists a closed definable set $W_P^S \subset \overline{\Delta_P^S}$ with ${\rm dim} (W_P^S) < k-1$ such that $\Delta_P^S \setminus W_P^S$ is smooth of dimension $k-1$ and such that the following function in $y \mapsto \alpha(\pi_P^S,y)$
is defined and constant on each connected component of $\Delta_P^S \setminus W_P^S$.
\end{lemma}
\proof Let 
$$W_P^S = \overline{ {\rm Sing}(\Delta_P^S) \cup \pi_P^S (T_P^S \cup T'^S_P) \cup \pi_P^S ({\rm Fr}(\Sigma_P^S))  }.$$
By the previous results, $W_P^S$ is a definable set of dimension strictly less than $k-1$. Moreover, $W_P^S$ is a closed set in $P$, hence $\Delta_P^S \setminus W_P^S$ is an open set in $\Delta_P^S$. 
The set $\Delta_P^S \setminus W_P^S$ is a smooth $(k-1)$-dimensional manifold for ${\rm Sing}(\Delta_P^S) \subset W_P^S$ and the function of the statement is well-defined because 
$\pi_P^S (T_P^S) \cup {\rm Sing}(\Delta_P^S ) \subset W_P^S$. 

Let $y$ be a point in $\Delta_P^S \setminus W_P^S$ and let $x$ be the point in $\Sigma_P^S \setminus (T_P^S \cup T'^S_P)$ such that $\pi_P^S(x)=y$. We can choose an open neighborhood $U_x \subset \Sigma_P^S \setminus (T_P^S \cup T'^S_P)$ such that ${\pi_P^S}_{\vert U_x}$ is a diffeomorphism and such that 
$$\alpha(\pi_P^S,y) =\alpha(\pi_P^S, \pi_P^S(x')),$$ for each $x'$ in $U_x$ by Lemma \ref{IndiceLocConstant}. Let $V \subsetneq U_x$ be a smaller open neighborhood of $x$ in $\Sigma_P^S \setminus (T_P^S \cup T'^S_P)$. Let $A$ be the following set:
$$A = \overline{\Sigma_P^S} \setminus V.$$
It is a compact subset of $\overline{\Sigma_P^S}$, hence $\pi_P^S(A)$ is compact in $\overline{\Delta_P^S}$. 
The point $y$ does not belong to $\overline{\pi_P^S(A)}$ because $y$ does not belong neither to $\pi_P^S ({\rm Fr}(\Sigma_P^S))$ nor to $\pi_P^S(T_P^S \cup T'^S_P)$. Hence there exists an open neighborhood $\mathcal{U}$ of $y$ in  $\overline{\Delta_P^S}$ which does not intersect $\pi_P^S(A)$. The function $y' \mapsto \alpha(\pi_P^S,y')$  is constant on $\mathcal{U}'=\mathcal{U} \cap (\Delta_P^S \setminus W_P^S)$. Note that $(\pi_P^S)^{-1}(\mathcal{U}')$ is included $U_x$.\endproof

Let $U$ be an open subset of $X$. Recall that $\Delta_P^{S \cap U}=\pi_P^S(\Sigma_P^{S\cap U}) = \pi_P^S(\Sigma_P^S \cap U)$.
\begin{lemma}\label{Open}
The set $\Delta_P^{S\cap U} \setminus W_P^S$ is an open subset of $\Delta_P^S \setminus W_P^S$.
\end{lemma}
\proof The inclusion is clear. Now let $y$ be a point in $\Delta_P^{S\cap U} \setminus W_P^S$. There exists a unique point $x$ in $(\Sigma_P^S \cap U) \setminus (T_P^S \cup T'^S_P)$ such that $\pi_P^S(x)=y$. We can find a small open neighborhood $U_x$ of $x$ in $\Sigma_P^S$ such that $U_x$ is included in 
$(\Sigma_P^S \cap U) \setminus (T_P^S \cup T'^S_P)$. As above, there is an open neighborhood $\mathcal{U}'$ of $y$ in $\Delta_P^S \setminus W_P^S$ such that $(\pi_P^S)^{-1}(\mathcal{U}')$ is included in $U_x$. \endproof

\begin{definition}\label{DefPolarLength1}
{\rm Let $U $ be an open subset of $X$ and let $S$ be a stratum of $X$. If $d_S < n$, for each $q \in \{0,\ldots,d_S\}$ and for $P \in G_n^{q+1}$ generic in the sense of Lemmas \ref{GenericProj}, \ref{GenericFolds}, \ref{GenericDoublePoints} and Proposition \ref{GenericProjLim}, we set
$$m_{S,q}(P,U) = \int_{\Delta_P^{S \cap U}} \alpha(\pi_P^S,y) dy .$$
We define
$$L_q(X,S,U)= \frac{\beta(1,n-q)}{\beta(q+1,n) g_n^{q+1} } \int_{G_n^{q+1}} m_{S,q} (P,U) dP.$$
For $d_S < q \le n$, we set $L_q(X,S,U)=0$. 

If $d_S=n$, we set $L_q(X,S,U)=0$ if $q < n$ and $L_n (X,S,U)= {\rm vol}(S \cap U)$.}
\end{definition}

We note that $m_{S,q}(P,U)$ is well-defined because the definable set $W_P^S$ has dimension at most $k-2$ and the index $\alpha(\pi_P^S,y)$ is defined for $y$ in $\Delta_P^{S\cap U} \setminus W_P^S \subset \Delta_P^S \setminus W_P^S$.

We are in position to define the polar lengths (or polar measures) of a definable stratified set. 
\begin{definition}\label{DefPolarLength2}
{\rm Let $X \subset \mathbb{R}^n$ be a compact definable set equipped with a finite definable stratification $X = \sqcup_{a \in A} S_a$. Let $U$ be an open subset of $X$. For each $q \in \{0,\ldots,n \}$, we set
$$L_q(X,U) = \sum_{a \in A} L_q(X,S_a,U).$$}
\end{definition}
First we remark that $L_0(X,U)= \Lambda_0(X,U)$ since
$$ L_0(X,U)= \frac{1}{g_n^1} \int_{G_n^1} \left( \sum_{x \in X} \frac{1}{2} \Big( {\rm ind}(v^*,X,x) + {\rm ind}(-v^*,X,x) \Big) \right) dL,$$
where $v$ is a unit vector such that $L=[ v ]$. This can be rewritten 
$$ L_0(X,U)= \frac{1}{s_{n-1}} \int_{S^{n-1}} \left( \sum_{x \in X} \frac{1}{2} \Big( {\rm ind}(v^*,X,x) + {\rm ind}(-v^*,X,x) \Big) \right) dv.$$
The right-hand side of this equality coincides with the right-hand side of the equality of Proposition \ref{ExchangeFormula}.

The polar lengths satisfy the linear kinematic formula.
\begin{proposition}\label{KinematicLength}
For $q \in \{0,\ldots,n\}$, we have
$$L_q(X,U)= {\rm cst} \int_{A_n^{n-q}} L_0 (X \cap L, X \cap L \cap U) dL.$$
\end{proposition}
\proof We assume first that $d_X<n$. If $q > d_X$ then the result is clear because almost all $L$ in $A_n^{n-q}$ does not intersect $X$.
Let us fix $q \le d_X$. We have
$$L_q(X,U) = \sum_{a \ \vert \ d_{S_a} \ge q} L_q (X,S_a, U).$$
Let us a fix a stratum $S$ such that $d_S \ge q$. We have
$$L_q(X,S,U)= \frac{\beta(1,n-q)}{\beta(q+1,n) g_n^{q+1} } \int_{G_n^{q+1}} \int_{\Delta_P^{S \cap U}}  \alpha(\pi_P^S,y) dy  dP.$$
It is clear that
$$ \int_{\Delta_P^{S \cap U}} \alpha(\pi_P^S,y) dy =  \int_{\Delta_P^{S \cap U} \setminus W_p^S} \alpha(\pi_P^S,y) dy.$$
Let us decompose $\Delta_P^{S} \setminus W_P^{S }$ into the finite union of its connected components, i.e. $\Delta_P^{S} \setminus W_P^{S}= \cup Y_j^{P,S}$. For each $j$, let us denote by $\lambda_j^{P,S}$ the value of 
$\alpha(\pi_P^S,y)$ on $Y_j^{P,S}$. We can write 
$$\int_{\Delta_P^{S \cap U}} \alpha(\pi_P^S,y) dy = \sum_j \lambda_j^{P,S} \cdot {\rm vol}(Y_j^{P,S} \cap \Delta_P^{S\cap U} ).$$
The Cauchy-Crofton formula \cite{Federer} in $P$  gives 
$${\rm vol}(Y_j^{P,S} \cap \Delta_P^{S\cap U}) = {\rm cst} \int_{A_P^1}  \# (Y_j^{P,S} \cap \Delta_P^{S\cap U} \cap l) dl,$$
and so 
$$L_q(X,S,U) = {\rm cst} \int_{G_n^{q+1}} \left[ \int_{A_P^1}  \sum_j \lambda_j^{P,S} \cdot  \# (Y_j^{P,S} \cap \Delta_P^{S\cap U} \cap l) dl \right] dP.$$
Let $l$ be a generic line in $A_P^1$ which intersects the $Y_j^{P,S}$'s transversally and does not meet $W_P^{S}$. Let us write $(\cup_j Y_j^{P,S}) \cap l = \{y_1,\ldots,y_t \}$. Then the affine space $L$, defined by $L = P^\perp \oplus l$ intersects $S$ transversally. 
To see this, let $x$ be a point in $S \cap L$. 
If $x \notin \cup_{k=1}^t ({\pi_P^{S}}_{\vert \Sigma_P^{S}} )^{-1}(y_k)$ then $P^\perp$ intersects $S$ transversally at $x$ and so $L$ intersects $S $ transversally at $x$. If $x \in \cup_{k=1}^t ({\pi_P^{S }}_{\vert \Sigma_P^{S }} )^{-1}(y_k)$ then $P^\perp \oplus l$ is exactly the affine space $Q_{u,x}$, where $u$ is a unit vector in the direction of $l$. Furthermore the set $\cup_{k=1}^t ({\pi_P^{S}}_{\vert \Sigma_P^{S}} )^{-1}(y_k)$ is exactly the set of critical points of $u^* : S  \cap L \rightarrow \mathbb{R}$. Let us denote by $\Gamma^S_U$ the set of these critical points lying in $S \cap U$. We get
$$\sum_j \lambda_j^{P,S} \cdot \#(Y_j^{P,S} \cap \Delta_P^{S \cap U} \cap l ) = \sum_{x \in \Gamma^S_U}  \alpha(\pi_P^S,\pi_P^S(x)) ,$$
and 
$$L_q(X,S,U) = {\rm cst} \int_{G_n^{q+1}} \int_{A_P^1} \sum_{ x \in \Gamma^S_U} \alpha(\pi_P^S,\pi_P^S(x)) dl dP.$$
Let $\mathcal{F}$ be the flag variety of pairs $(P,l)$, $P \in G_n^{q+1}$ and $l \in A_P^1$. The mapping $(P,l) \mapsto (L,l)$, where $L=P^\perp \oplus l$, enables us to identify $\mathcal{F}$ with the flag variety of pairs $(L,l)$, where $L \in A_n^{n-q}$ and $l \in G_L^1$. Therefore we can write 
$$L_q(X,S,U)= {\rm cst} \int_{A_n^{n-q}} \left( \int_{G_L^1} \sum_{x \in \Gamma^S_U}  \alpha(\pi_P^S,\pi_P^S(x)) dl \right) dL,$$
where $u$ is a unit vector vector in $l$ and $P=L^\perp \oplus l$. But for almost all $u$ unit vector in $L$, $u^*_{\vert X \cap L}$ is a Morse function (see \cite{BroeckerKuppe}, Lemma 3.5), and by Lemma \ref{MorseIndex2}, 
$$\alpha(\pi_P^S,\pi_P^S(x))= \frac{1}{2} \Big( {\rm ind}(u^*,X \cap L,x)  + {\rm ind}(-u^*,X \cap L,x) \Big).$$
Hence by the definition of $L_0$, we can conclude that 
$$L_q(X,S,U) ={\rm cst} \int_{A_n^{n-q}} L_0(X\cap L, S \cap L,  X \cap L \cap U) dL,$$
and so
$$L_q (X,U)= \sum_{a \ \vert \ d_{S_a} \ge q} {\rm cst} \int_{A_n^{n-q}} L_0(X\cap L, S_a \cap L, X \cap L \cap U) dL=$$
$${\rm cst} \sum_a  \int_{A_n^{n-q}} L_0(X\cap L, S_a \cap L, X\cap L \cap U) dL,$$
because if $d_{S_a} < q$ then generically $L$ does not meet $S_a$. Finally we can write 
$$L_q (X,U) = {\rm cst} \int_{A_n^{n-q}} L_0(X\cap L, X \cap L \cap U) dL.$$
If $d_X=n$ then for $q<n$
$$\displaylines{
\quad L_q(X,U)= \sum_{a \ \vert \ d_{S_a} < n}    L_q(X,S_a,U) = \hfill \cr
\hfill  {\rm cst}   \sum_{a \ \vert \ d_{S_a} < n}\int_{A_n^{n-q}} L_0(X\cap L, S_a \cap L, X \cap L \cap U) dL. \quad \cr
}$$
But if $S$ is a stratum of dimension $n$ then $L_0(X \cap L, S \cap L, X \cap L \cap U)=0$ because $S \cap L \cap U$ is an open subset of the affine space $L$, on which any $u^*$ has no critical point. Furthermore
$$L_n(X,U) = \sum_{a \ \vert \ d_{S_a} = n} L_n(X,S,U) = \sum_{a \ \vert \ d_{S_a} =n} {\rm vol}(S_a \cap U) =$$ $$  \sum_{a \ \vert \ d_{S_a} =n} \int_{S_a \cap U} dx = \sum_{a \ \vert \ d_{S_a} =n} \int_{\mathbb{R}^n} L_0 (X \cap \{x\}, S_a \cap \{x\}, X \cap \{x\} \cap U) dx .$$

\endproof

\begin{theorem}\label{Curvature=Length}
For $q \in \{0,\ldots,n\}$, $L_q (X,U)= \Lambda_q(X,U)$. The polar lengths $L_q(X,U)$ do not depend on the stratification of $X$.
\end{theorem}
\proof First let us recall that $L_0(X,U)= \Lambda_0(X,U)$. By the linear kinematic formula for $\Lambda_q$ (see Proposition \ref{GlobalKinematicFormula}) and the previous proposition, we see that $L_q$ and $\Lambda_q$ are proportional. Let us prove that they coincide on $q$-dimensional smooth sets. It is true by definition for $n$-dimensional smooth sets.  Let $S$ be a smooth compact definable manifold of dimension $q < n$ and let $U$ be an open subset of $S$. By the Cauchy-Crofton formula \cite{Federer}, we have
$${\rm vol}(S \cap U) = \frac{1}{\beta(q,n) g_n^q} \int_{G_n^q}  \left[ \int_H \left( \# \pi_H^{S \cap U} \right)^{-1}(z) dz \right] dH,$$
and
$${\rm vol}(S \cap U) = \frac{1}{\beta(q,n) g_n^q} \int_{G_n^q}  \left[ \int_H   \frac{1}{g_{n-q}^1} \int_{G_{H^\perp}^1} \left( \# \pi_H^{S \cap U} \right)^{-1}(z) dl dz \right] dH.$$
Let $P \in G_n^{q+1}$ such that $H \subset P$ and let $\pi_H^P$ be the orthogonal projection on $H$ restricted to $P$. We have $\pi_H^{S \cap U} = \pi_H^P \circ \pi_P^{S \cap U}$ and so
$$\#(\pi_H^{S \cap U})^{-1} (z) = \# \left\{(\pi_P^{S \cap U})^{-1} \left[ (\pi_H^P)^{-1}(z) \right] \right\},$$
i.e.,
$$\#(\pi_H^{S \cap U})^{-1} (z) = \sum_{y \in (\pi_H^P)^{-1}(z) \cap \Delta_P^{S \cap U}} \# (\pi_P^{S \cap U})^{-1} (y).$$
Writing $P= H \oplus l$, where $l \in G_{H^\perp}^1$, we can identify the flag variety of pairs $(P,H)$, $P \in G_n^{q+1}$ and $H \in G_P^q$, with the flag variety $(H,l)$, $H \in G_n^q$ and $l \in G_{H^\perp}^1$.   Therefore, we can write
$${\rm vol}(S \cap U) = \frac{1}{\beta(q,n) g_n^q g_{n-q}^1} \int_{G_n^{q+1}} \int_{G_P^q} \int_H 
\sum_{y \in (\pi_H^P)^{-1}(z) \cap \Delta_P^{S \cap U}} \# (\pi_P^{S \cap U})^{-1}(y) dz dH dP.$$
Now we know that for $P$ generic in $G_n^{q+1}$, there is a definable set $W_P^S$ of dimension strictly less than $q$ such that $\# (\pi_P^S)^{-1}(y)$ is equal to $1$ on $\Delta_P^S \setminus W_P^S$. Therefore we have
$$\sum_{y \in (\pi_H^P)^{-1}(z)\cap \Delta_P^{S \cap U} } \# (\pi_P^{S \cap U})^{-1}(y) = \# \left\{ (\pi_H^P)^{-1}(z) \cap \Delta_P^{S \cap U} \right\}.$$
Applying the Cauchy-Crofton formula to compute vol$(\Delta_P^{S \cap U})$, we find that
$${\rm vol}(S \cap U) = \frac{\beta(q,q+1) g_{q+1}^q }{\beta(q,n) g_n^1 g_{n-q}^1} \int_{G_n^{q+1}} {\rm vol}(\Delta_P^{S \cap U}) dP,$$
which we can rewrite
$${\rm vol}(S \cap U) = \frac{\beta(1,n-q)}{\beta(q+1,n) g_n^{q+1}}  \int_{G_n^{q+1}} \int_{\Delta_P^{S \cap U}} \alpha(\pi_P^S,y) dy dP= L_q(S,U),$$
because here $\alpha(\pi_P^S,y)=1$.  \endproof
Let us remark that this last equality can be stated in a nicer way. Let $X \subset \mathbb{R}^n$ be a compact definable set and let $U \subset X$ be an open subset. If $d_X \le n-1$ then
$${\rm vol}(X \cap U)= \frac{\beta(1,n-d_X)}{\beta(d_X+1,n) g_n^{d_X+1}}  \int_{G_n^{d_X+1}}   {\rm vol}(\pi_P^X(X \cap U)) dP.$$
Of course, this equality is trivial if $d_X=n-1$.

\section{Limits of Lipschitz-Killing curvatures and density of polar images}
In this section, we give formulas which can be viewed as infinitesimal versions of the formulas established in the previous section. 

Let $(X,0) \subset (\mathbb{R}^n,0)$ be the germ at the origin of a closed definable set. We assume that $X$ is equipped with a finite definable stratification $\{ S_a \}_{a = 0}^t$, where $S_0$ is the stratum that contains $0$ and $0$ belongs to the frontier of each $S_a$, $a=1,\ldots,t$. 

As in the previous section, for every stratum $S$ of $X$, $\Sigma_P^S$ will denote the polar variety of $\pi_P^S$, where $P \in G_n^k$, $k=1,\ldots,n$. As explained in \cite{ComteMerle}, the projection of a definable germ does not define necessarily a germ. But, generically, it defines a definable germ thanks to the following proposition.

\begin{proposition}\label{local1}
Let $k \in \{1,\ldots,n\}$. For almost all $P \in G_n^k$, there exists an open neighborhood $U$ of $0$ in $\mathbb{R}^n$ such that
$$(U \cap \overline{\Sigma_P^S} \cap P^\perp) \setminus \{0\} = \emptyset.$$
\end{proposition}
\proof See \cite{ComteMerle}, Proposition 2.2. \endproof

\begin{corollary}\label{local2}
Let $(X,0) \subset (\mathbb{R}^n,0)$ be the germ of a closed definable set. Let $k \in \{1,\ldots,n\}$. Then for almost all $P \in G_n^k$, there exits $\epsilon_P >0$ such that for all $\epsilon \in ]0,\epsilon_P]$, there exists $\eta_\epsilon >0$ such that
$$B^P_{\eta_\epsilon} \cap \pi_P^X (X \cap B_\epsilon) = B^P_{\eta_\epsilon} \cap \pi_P^X (X \cap B_{\epsilon_P}),$$
where $B^P_{\eta_\epsilon}$ is the ball of radius $\eta_\epsilon$ centered at the origin in $P$.
In other words, a generic projection of the germ $(X,0)$ defines a germ at the origin in $P$.
\end{corollary}
\proof See \cite{ComteMerle}, Proposition 2.3. \endproof

\begin{corollary}\label{local3}
Let $S$ be a stratum of $X$. Then for almost all $P \in G_n^k$, there exists $\epsilon_P >0$ such that for all $\epsilon \in ]0,\epsilon_P]$, there exists $\eta_\epsilon$ such that
$$B^P_{\eta_\epsilon} \cap \pi_P^X (\overline{\Sigma_P^S}  \cap B_\epsilon) = B^P_{\eta_\epsilon} \cap \pi_P^X (\overline{\Sigma_P^S}  \cap B_{\epsilon_P}),$$
and
$$B^P_{\eta_\epsilon} \cap \pi_P^X ({\rm Fr}(\overline{\Sigma_P^S}  \cap B_\epsilon)) = B^P_{\eta_\epsilon} \cap \pi_P^X ({\rm Fr}(\overline{\Sigma_P^S}  \cap B_{\epsilon_P})).$$
\end{corollary}
\proof See \cite{ComteMerle}, Proposition 2.4. \endproof

For each stratum $S$ of $X$ and for $P$ generic in $G_n^k$, $k=1,\ldots,d_S+1$, we recall that $T_P^S$ is the following closed definable subset of $\Sigma_P^S$ of dimension at most $k-2$: 
$$ T_P^S =Z_P^S \bigcup \overline{D({\pi_P^S}_{\vert \Sigma_P^{S,0}})}  \bigcup \overline{( \cup_i \Sigma_{P,\lim}^{S_i}) \cap S},$$ where the strata $S_i$ contain $S$ in their frontier, and that $T'^S_P$ is the closed definable subset of $\Sigma_P^S$ of dimension at most $k-2$ defined in Lemma \ref{IndiceLocConstant}.
As above, we have:
\begin{corollary}\label{local4}
Let $S$ be a stratum of $X$ and let $k \in \{1,\ldots,d_S+1\}$. Then for almost all $P \in G_n^k$, there exists $\epsilon_P >0$ such that for all $\epsilon \in ]0,\epsilon_P]$, there exists $\eta_\epsilon$ such that
$$B^P_{\eta_\epsilon} \cap \pi_P^X (\overline{T_P^S \cup T'^S_P}  \cap B_\epsilon^n) = B^P_{\eta_\epsilon} \cap \pi_P^X (\overline{T_P^S \cup T'^S_P}   \cap B_{\epsilon_P}^n).$$
\end{corollary}

For $\epsilon >0$, let us denote by $\Delta_P^{S,\epsilon}$ the definable set $\pi_P^S(\overline{\Sigma_P^S} \cap B_\epsilon^n)$. By the results of Section 2, generically, $\Delta_P^{S,\epsilon}$ is a definable subset of $P$ of dimension $k-1$. 
Let $W_P^{S,\epsilon}$ be the following closed definable set:
$$W_P^{S,\epsilon} = \overline{ {\rm Sing}(\Delta_P^{S,\epsilon}) \cup \pi_P^S ((T_P^S\cup T'^S_P) \cap B_\epsilon^n) \cup \pi_P^S ({\rm Fr}(\overline{\Sigma_P^S} \cap B_\epsilon^n ))  }.$$
As in Section 3, $W_P^{S,\epsilon}$ is a closed definable set included in $\Delta_P^{S,\epsilon}$ of dimension $k-2$.
Let us consider the connected components of $\Delta_P^{S,\epsilon} \setminus W_P^{S,\epsilon}$ that contain $0$ in their closure. We denote them by $Y_1^{P,S,\epsilon},\ldots,Y_{r_{P,S}}^{P,S,\epsilon}$. As in Section 2, on each $Y_j^{P,S,\epsilon}$ the function $ y \mapsto \alpha(\pi_P^S,y)$ is constant.
Furthermore, since by Corollaries \ref{local3} and \ref{local4} the sets $\Delta_P^{S,\epsilon}$ and $W_P^{S,\epsilon}$ define germs, the sets $Y_1^{P,S,\epsilon},\ldots,Y_{r_P}^{P,S,\epsilon}$ define germs as well. 

From now on, we will denote by $\Delta_P^S$ the germ defined by the $\Delta_P^{S,\epsilon}$'s and for $j \in \{1,\ldots,r_{P,S} \}$, by $Y_j^{P,S}$ the germ defined by the $Y_j^{P,S,\epsilon}$'s. To each $Y_j^{P,S}$, we can associate the following integer 
$$\lambda_j^{P,S}= \lim_{\epsilon \rightarrow 0} \lim_{y \rightarrow 0 \atop y \in Y_j^{P,S}}  \alpha(\pi_P^S,y).$$
This integer does not depend on $\epsilon$ nor on $y$, provided they are small enough. We are now in position to define the localizations of the polar lengths.
\begin{definition}\label{LocLength1}
{\rm Let $S$ be a stratum of $X$ such that $0 \in \overline{S}$. If $d_S <n$ then for each $q \in \{0,\ldots,d_S\}$, we set 
$$L_q^{\rm loc}(X,S,0)= \frac{1}{g_n^{q+1}} \int_{G_n^{q+1}} \left[ \sum_{j=1}^{r_{P,S}} \lambda_j^{P,S} \cdot \Theta_q(Y_j^{P,S},0) \right] dP.$$
For $q > d_S$, we set $L_q^{\rm loc}(X,S,0)=0$.

If $d_S=n$ then we set $L_q^{\rm loc}(X,S,0)=0$ for $q<n$ and $L_n^{\rm loc}(X,S,0) = \Theta_n(S,0)$.}
\end{definition}

\begin{definition}\label{LocLength2}
{\rm Let $(X,0)$ be the germ of a closed definable set and let $\{S_a\}_{a=0}^t $ be a definable stratification of $X$ such that $0 \in \overline{S_a}$ for $\alpha \in \{0,\ldots,t\}$. For each $q \in \{0,\ldots,n\}$, we set 
$$L_q^{\rm loc} (X,0)= \sum_{a=0}^t L_q^{\rm loc} (X,S_a,0).$$}
\end{definition}
The following theorem relates the limits $\lim_{\epsilon \to 0} \frac{\Lambda_k(X,X \cap B_\epsilon ^n)}{b_k \epsilon^k}$ to the polar lengths. It can be viewed as an infinitesimal version of Theorem \ref{Curvature=Length}.
\begin{theorem}\label{Curvature=LengthLocal}
Let $(X,0) \subset (\mathbb{R}^n,0)$ be the germ of a closed definable set. For $k \in \{0,\ldots,n\}$, we have
$$\lim_{\epsilon \rightarrow 0} \frac{\Lambda_k(X,X\cap B_\epsilon)}{b_k \epsilon^k} = L_k^{\rm loc}(X,0).$$
\end{theorem}
\proof  Let us assume first that $d_X <n$. If $k > d_X$ then the result is obvious by the definitions of the $\Lambda_k$'s and $L_k^{\rm loc}$'s. Let us fix $k \le d_X$. In this case,
$$L_k^{\rm loc}(X,0)= \sum_{a \ \vert \ d_{S_a} \ge k } L_k^{\rm loc}(X,S_a,0).$$
Let us fix a stratum $S$ such that $d_S \ge k$. We know that 
$$L_k^{\rm loc}(X,S,0) = \frac{1}{g_n^{k+1}} \int_{G_n^{k+1}} \left[ \sum_{j=1}^{r_{P,S}} \lambda_j^{P,S} \cdot \Theta_k (Y_j^{P,S},0) \right] dP.$$
We apply the Cauchy-Crofton formula for the density due to Comte \cite{Comte}. Actually we will use here a bit different version of this formula, namely we will apply Theorem \ref{LocalKinematicFormula}. 

For $P \in G_n^{k+1}$, let $L \in G_P^1$ and $L^{\perp P}$ be the orthogonal complement of $L$ in $P$. For each $v$ in $S_{L^{\perp P}}$ and each $\delta >0$, let $L_{\delta,v}$ be the affine line $L+\delta v$. We set 
$$\beta_{0,j}(L,v) = \lim_{\epsilon \rightarrow 0}  \lim_{\delta \rightarrow 0}  \# L_{\delta,v} \cap Y_j^{P,S} \cap B_\epsilon^P.$$
By this theorem, we know that
$$\Theta_k(Y_j^{P,S},0) = \frac{1}{g_{k+1}^1} \int_{G_P^1} \left( \frac{1}{s_{k-1}} \int_{S_{L^\perp P}} \beta_{0,j}(L,v) dv \right) dL.$$
In fact, Theorem \ref{LocalKinematicFormula} is proved  in \cite{DutertreProcTrotman} for germs of closed sets. But, it holds also here because $\Theta_k(Y_j^{P,S},0) = \Theta_k(\overline{Y_j^{P,S}},0)$ and generically $\# L_{\delta,v} \cap {\rm Fr}(\overline{Y_j^{P,S}}) =0$. Then we find
$$L_k^{\rm loc}(X,S,0) = \frac{1}{g_n^{k+1}} \int_{G_n^{k+1}} \frac{1}{g_{k+1}^1} \int_{G_P^1} \frac{1}{s_{k-1}} \int_{S_{L^{\perp P}}} I^S(P,L,v) dv dL dP,$$
where $I^S(P,L,v) =\sum_{j=1}^{r_{P,S}} \lambda_j^{P,S} \beta_{0,j} (L,v)$. Since the flag variety of pairs $(P,L)$, $P \in G_n^{k+1}$ and $L \in G_P^1$, is isomorphic to the flag variety of pairs $(H,L)$, $H \in G_n^{n-k}$ and $L \in G_H^1$, by the mapping $(P,L) \mapsto (H,L)$ where $H = P^\perp \oplus L$, we can write
$$L_k^{\rm loc}(X,S,0) =\frac{1}{g_n^{k+1} g_{k+1}^1} \int_{G_n^{n-k}} \frac{1}{s_{k-1}} \int_{S_{H^\perp}} \left( \int_{G_H^1} I^S(H,L,v) dH \right) dv dP,$$
because $H^\perp= L^{\perp P}$, and where $I^S(H,L,v)=I^S(P,L,v)$. As in the global case, we see that the integer $I^S(H,L,v)$ admits the following description:
$$ I^S(H,L,v) =\lim_{\epsilon \to 0} \lim_{\delta \to 0} \sum_{x \in \Gamma_{\delta,\epsilon}^S} \alpha(\pi_P^S, \pi_P^S(x)),$$
where $0 < \vert \delta \vert \ll \epsilon \ll 1$ and $\Gamma_{\delta,\epsilon}^S$ is the set of critical points of $u^* : S \cap H_{\delta,v} \cap \mathring{B_\epsilon} \rightarrow \mathbb{R}$, $u$ being a unit vector of the line $L$.

Let us choose a sequence $(\epsilon_m)_{m \in \mathbb{N}}$ such that $\lim_{m \rightarrow + \infty} \epsilon_m =0$ and for each $m \in \mathbb{N}$, a sequence and $(\delta_p^m)_{p \in \mathbb{N}}$ such that $\lim_{p \rightarrow + \infty} \delta_p^m =0$.
By \cite{BroeckerKuppe} Lemma 3.5, for almost all $L$ in $G_H^1$, $u^*_{\vert X \cap H_{\delta_p^m,v} \cap B_{\epsilon_m}}$ is a Morse function. Since a countable union of sets of measure zero has measure zero, for almost all $L$ in $G_H^1$, the function $u^*_{\vert X \cap H_{\delta_p^m,v} \cap B_{\epsilon_m}}$ is a Morse function for every $(m,p) \in \mathbb{N}^2$.
Hence by Lemma \ref{MorseIndex2}, for every $x \in \Gamma^S_{\delta_p^m,\epsilon^m}$, we have
$$\alpha(\pi_P^S, \pi_P^S(x))= \frac{1}{2} \left( {\rm ind}(u^*,X \cap H_{\delta_p^m,v},x) + {\rm ind}(-u^*,X \cap H_{\delta_p^m,v},x) \right).$$
Therefore we find that
$$\displaylines{
\quad \int_{G_H^1} I^S(H,L,v) dL = \hfill \cr
  \int_{G_H^1} \lim_{m \rightarrow + \infty}  \lim_{p \to \infty} \sum_{x \in \Lambda^S_{\delta_p^m,\epsilon_m}} \frac{1}{2}   \left( {\rm ind}(u^*,X \cap H_{\delta_p^m,v},x) + {\rm ind}(-u^*,X \cap H_{\delta_p^m,v},x) \right) dL= \cr$$
\hfill g_{n-k}^1 \lim_{m \rightarrow + \infty}  \lim_{p \to \infty} \Lambda_0(X \cap H_{\delta_p^m,v} , S \cap H_{\delta_p^m,v} \cap B_{\epsilon_m} ). \quad \cr
}$$
Summing over all the strata $S_a$ such that $d_{S_a} \ge k$ leads to
$$L_k^{\rm loc}(X,0) = \frac{g_{n-k}^1}{g_n^{k+1} g_{k+1}^1 s_{k-1}}  \int_{S_{H^\perp}}  \lim_{m \rightarrow + \infty}  \lim_{p \to \infty}  \Lambda_0(X \cap H_{\delta_p^m,v} , X \cap H_{\delta_p^m,v} \cap B_{\epsilon_m} ) dv dH=$$
$$ \frac{1}{g_n^{n-k}} \int_{G_n^{n-k}} \frac{1}{s_{k-1}} \int_{S_{H^\perp}}  \lim_{\epsilon \rightarrow 0} \lim_{\delta \rightarrow 0} \Lambda_0(X \cap H_{\delta,v} , X \cap H_{\delta,v} \cap B_{\epsilon} ) dv dH.$$
Theorem \ref{LocalKinematicFormula} implies that 
$$L_k^{\rm loc}(X,0) = \lim_{\epsilon \rightarrow 0} \frac{\Lambda_k(X,X\cap B_\epsilon)}{b_k \epsilon_k}.$$
If $d_X=n$ and $k<n$ then the proof works in the same way because the strata of dimension $n$ have no contribution. If $k=n$ then the formula is obvious because both expressions are equal to the density of $(X,0)$.
\endproof

Applying Theorem \ref{CurvAndPolar}, we obtain the following relation between the polar invariants and the $L_k^{\rm loc}(X,0)$'s. 
\begin{corollary}\label{PolarLocalAndPolarInv}
Let $(X,0)$ be the germ of a closed definable set. For $k \in \{0,\ldots,n-1\}$, we have
$$\sigma_k(X,0) -\sigma_{k+1}(X,0) = L_k^{\rm loc}(X,0).$$
Furthermore, we have
$$\sigma_n(X,0) =L_n^{\rm loc}(X,0).$$
\end{corollary}
$\hfill \Box$

Let us refine the above results. By Definitions \ref{LocLength1} and \ref{LocLength2}, we know that $L_k^{\rm loc}(X,0)=0$ for $k >d_X$. Furthermore, it is explained in \cite{ComteMerle} that for $k \in \{1,\ldots,d_{S_0}\}$ and $P$ generic in $G_n^k$, the projection $\pi_P^X$ is a submersion in a neighborhood of $0$. Hence, generically the polar varieties and the polar images are empty near the origin and so $L_k^{\rm loc}(X,0)=0$ for $k\in \{0,\ldots,d_{S_0}-1\}$. 
Let us focus now on the localized polar length $L_{d_{S_0}}^{\rm loc}(X,0)$. If $d_S > d_{S_0}$, then generically ${\rm Fr}(\Delta_P^S)$ has dimension $d_{S_0}-1$ and $0$ is not in the closure of $\Delta_P^S$, so $L_{d_{S_0}}^{\rm loc}(X,S,0)=0$. Combining this with Definition \ref{LocLength2} and Corollary \ref{PolarLocalAndPolarInv}, we find that
$$L_{d_{S_0}}^{\rm loc}(X,S_0,0) =1-\sigma_{d_{S_0}+1}(X,0),$$
because $\sigma_{d_{S_0}}(X,0)=1$ (see \cite{ComteMerle}, Remarque 2.9).

If $d_{S_0}<k\le d$ then $S_0$ has no contribution in the computation of $L_k^{\rm loc}(X,0)$ and so $L_k^{\rm loc}(X,0)= \sum_{a=1}^t L_k^{\rm loc} (X,S_a,0)$.  Similarly the strata of dimension strictly less than $d_X$ have no contribution in the computation of $L_{d_X}^{\rm loc}(X,0)=\Theta_{d_X}(X,0)$. Moreover, if $S$ is a stratum of dimension $d_X$, the indices $\lambda_j^{P,S}$ that appear in the definition of $L_{d_X}^{\rm loc}(X,0)$ are generically equal to $1$. Therefore we can restate Corollary \ref{PolarLocalAndPolarInv}.
\begin{corollary}\label{PolarLocalAndPolarInv2}
Let $(X,0)$ be the germ of a closed definable set, equipped with a finite definable Whitney stratification $\{ S_a \}_{a = 0}^t$, where $S_0$ is the stratum that contains $0$ and $0$ belongs to the frontier of each $S_a$, $a=1,\ldots,t$. We have 
$$ L_{d_{S_0}}^{\rm loc}(X,0) =L_{d_{S_0}}^{\rm loc}(X,S_0,0) =1-\sigma_{d_{S_0}+1}(X,0).$$
For $d_{S_0} < k <d_X$, we have
$$\sigma_k(X,0) -\sigma_{k+1}(X,0) = \sum_{a=1}^t L_{k}^{\rm loc}(X,S_a,0).$$
Furthermore if $d_X <n$, then we have
$$\Theta_{d_X}(X,0)=L_{d_{X}}^{\rm loc}(X,0)  = \frac{1}{g_n^{{d_X}+1}} \int_{G_n^{d_X+1}} \Theta_{d_X}(\pi_P^X(X),0) dP.$$
\end{corollary}
Of course, this  last equality  is trivial if $d_X=n-1$.

In \cite{LeTeissierAnnals}, Theorem 6.1.9 and  \cite{LeTeissierArcata}, Theorem 4.1.1, L\^e and Teissier proved a relation between the vanishing Euler characteristics and the multiplicities of the polar varieties of a complex analytic germ. 
More precisely, if $(X,0) \subset (\mathbb{C}^n,0)$ is a complex analytic germ, they showed that the difference of two consecutive vanishing Euler characteristics of $X$ at $0$ is equal to a weighted sum of polar multiplicities, the sum been taken over the strata that contain $0$ in their frontier. Let us describe the weights that appear in this sum. If $S_0$ denotes the stratum that contains $0$ and $S_a$ is a stratum such that $0 \in {\rm Fr}(S_a)$, then the weight attached to $S_a$ is 
$$(-1)^{d_{S_a}-d_{S_0}-1} \cdot (1-\chi ({\rm Lk}^{\mathbb{C}} (X,S_a) ),$$
 where ${\rm Lk}^{\mathbb{C}}(X,S_a)$ is the complex link of the stratum $S_a$ in $X$ (see for instance \cite{GoreskyMacPherson}, II 2.2, for the definition of the complex link). The first term of this product is actually a complex tangential Morse index and the second term is a complex normal Morse index. In the equalities
$$\sigma_k(X,0) -\sigma_{k+1}(X,0) = \sum_{a=1}^t L_{k}^{\rm loc}(X,S_a,0), \  d_{S_0} < k <d_X,$$
the weights $\lambda_j^{P,S}$ are half-sums of two stratified Morse indices, which are also decomposed into the product of a tangential Morse index and a normal Morse index.
Therefore, and since the polar invariants are  real versions of the vanishing Euler characteristics, we can view the equalities 
$$\sigma_k(X,0) -\sigma_{k+1}(X,0) = \sum_{a=1}^t L_{k}^{\rm loc}(X,S_a,0), \  d_{S_0} < k <d_X,$$
as  real versions of L\^e and Teissier's results.

In \cite{ComteMerle}, Theorems 4.9 and 4.10, Comte and Merle showed that the polar invariants were continuous along the strata of a $(w)$-stratification of a closed subanalytic set. Recently, Nguyen and Valette \cite{NguyenValette} extended this result for a  $(b)$-stratification and for a closed set $X$, definable in a polynomially bounded structure. They even proved that the polar invariants were locally lipschitz along the strata of a $(w)$-stratification of $X$.  Note that these results are not true if the structure is not polynomially bounded. A counter-example to the continuity of the density along the strata of a $(b)$-stratified set, definable in a non-polynomially bounded structure, was found recently by Trotman and Valette (see \cite{Nguyen}). 
As a straightforward consequence of these results and Corollary \ref{PolarLocalAndPolarInv}, we see that  that localized polar lengths are continuous along the strata of a $(b)$-stratification of a closed definable set, and locally lipschitz if the stratification is $(w)$-regular. Let us examine these last results in the case of two strata.

So let $S$ be stratum  that contains $S_0$ in its frontier and such that the depth of $S_0$ in $\overline{S}$ is one. This means that there is no stratum $S'$ such that $S_0 \subset {\rm Fr}(S')$ and $S' \subset {\rm Fr}(S)$. Adapting the notations introduced before  Definition \ref{LocLength1} in an obvious way, we have that for $k \in \{d_{S_0}+1,\ldots,d_S\}$ and $z \in S_0$,
$$L_k^{\rm loc}(\overline{S},z)=\frac{1}{g_n^{k+1}} \int_{G_n^{k+1}} \sum_{j=1}^{r_{P,S}(z)} \lambda_j^{P,S}(z) \Theta \Big( Y_j^{P,S}(z), z \Big) dP,$$
where the $Y_j^{P,S}(z)$'s are defined in the same way as the $Y_j^{P,S}$ were defined, replacing the ball $B_\epsilon$ with the ball $B_\epsilon(z)$.  The index $\lambda_j^{P,S}(z)$ was defined like that:
$$\lambda_j^{P,S}(z)= \lim_{\epsilon \rightarrow 0} \lim_{y \rightarrow z \atop y \in Y_j^{P,S,\epsilon}(z)}  \alpha(\pi_P^S,y),$$
with 
$$\alpha(\pi_P^S,y)=\frac{1}{2} \Big( {\rm ind}(\nu^*, \overline{S} \cap Q_{\nu,x},x) +{\rm ind}(-\nu^*, \overline{S} \cap Q_{\nu,x},x) \Big),$$
where $\nu$ is a unit normal vector to $T_y \Delta_P^S$ in $P$ and $y= \pi_P^S(x)$. But since $S$ is a top stratum in $\overline{S}$, the indices $ {\rm ind}(\nu^*, \overline{S} \cap Q_{\nu,x},x)$ and $ {\rm ind}(-\nu^*, \overline{S} \cap Q_{\nu,x},x)$ are actually tangential Morse indices. They are equal if ${\rm dim}(S \cap Q_{\nu,x})$ is even and opposite if ${\rm dim}(S \cap Q_{\nu,x})$ is odd. The dimension of $S \cap Q_{\nu,x}$ is $d_S-k$, so $L_k^{\rm loc}(\overline{S},z)=0$ if $d_S-k$ is odd and if $d_S-k$ is even,
$$L_k^{\rm loc}(\overline{S},z)=\frac{1}{g_n^{k+1}} \int_{G_n^{k+1}} \sum_{j=1}^{r_{P,S}(z)} \lambda_j^{P,S}(z) \Theta \Big( Y_j^{P,S}(z), z \Big) dP,$$
where $\lambda_j^{P,S}(z) \in \{-1,+1\}$. Hence in this situation, the localized polar lengths $L_k^{\rm loc}(\overline{S},-)$ are somehow real counterparts of the complex polar multiplicities. 
\begin{proposition}
Assume that the depth of $S_0$ in $\overline{S}$ is one. In this case, the localized polar lengths $L_k^{\rm loc}(\overline{S},-)$, $k \in \{d_{S_0}+1,\ldots, d_S\}$ and $d_S-k$ even, are continuous along $S_0$ if the stratification is $(b)$-regular and locally lipschitz if the stratification is $(w)$-regular.
\end{proposition}
\proof By Corollary \ref{PolarLocalAndPolarInv2}, for $d_{S_0} < k  <d_S$ we have
$$\sigma_k(\overline{S},z) -\sigma_{k+1}(\overline{S},z) =L_k^{\rm loc}(\overline{S},z)=L_k^{\rm loc}(\overline{S},S,z),$$
and $\Theta_{d_S}(\overline{S},z)= L_{d_S}^{\rm loc}(\overline{S},S,z)=L_{d_S}^{\rm loc}(\overline{S},z).$
Then we apply the results of \cite{NguyenValette}. \endproof
It would be interesting to see if this result still holds without the assumption on the depth of $S_0$ in $\overline{S}$.

\end{document}